\documentclass[12pt]{amsart}
\usepackage{amscd,amssymb,hyperref}
\usepackage[PostScript=dvips]{diagrams}
\usepackage{graphicx}
\usepackage{graphics}
\usepackage{epsfig}
\usepackage{picinpar}
\usepackage{wrapfig}
\usepackage{amsmath}
\usepackage{color}
\usepackage{pstricks,pst-node,pst-text,pst-3d}
\usepackage{pst-plot}
\usepackage{verbatim}

\newtheorem{thm}{Theorem}[section]
\newtheorem{lem}[thm]{Lemma}
\newtheorem{cor}[thm]{Corollary}
\newtheorem{prop}[thm]{Proposition}

\newtheorem{dfn}[thm]{Definition}

\newtheorem{rem}[thm]{Remark}

\def\square{\vbox{
      \hrule height 0.4pt
      \hbox{\vrule width 0.4pt height 5.5pt \kern 5.5pt \vrule width 0.4pt}
      \hrule height 0.4pt}}

\def\id{\mathrm{id}}

\def\Ker{\mathrm{K er}}
\def\ch\mathrm{c h}

\def\RP{\mathbb{R}\mathrm{P}}

\newcommand{\Z}{\mathbb{Z}}

\newcommand{\Brun}{\mathrm{Brun}}

\let\la=\langle
\let\ra=\rangle

\numberwithin{equation}{section}

\newcommand{\auths}[1]{\textrm{#1},}
\newcommand{\artTitle}[1]{\textsl{#1},}
\newcommand{\jTitle}[1]{\textrm{#1}}
\newcommand{\Vol}[1]{\textbf{#1}}
\newcommand{\Year}[1]{\textrm{(#1)}}
\newcommand{\Pages}[1]{\textrm{#1}}

\title{Brunnian Braids on Surfaces}

\author{V. G. Bardakov}
\address{Sobolev Institute of Mathematics, Novosibirsk 630090, Russia}
\email{bardakov@math.nsc.ru}

\author{R. Mikhailov}
\address{Steklov Mathematical Institute, Gubkina 8, 119991 Moscow, Russia}
\email{rmikhailov@mail.ru}

\author{V. V. Vershinin}
\address{D\'epartement des Sciences Math\'ematiques,
                                     Universit\'e Montpellier II,
Place Eug\`ene Bataillon,
34095 Montpellier cedex 5, France}
\email{ vershini@math.univ-montp2.fr}
\address{Sobolev Institute of Mathematics, Novosibirsk 630090,
Russia }
\email{ versh@math.nsc.ru}

\author{J. Wu}
\address{Department of Mathematics, National University of Singapore, 2 Science Drive 2
Singapore 117542} \email{matwuj@nus.edu.sg}
\urladdr{www.math.nus.edu.sg/\~{}matwujie}



\begin{document}

\begin{abstract}
We determine a set of generators for the Brunnian braids on a general surface $M$ for $M\not=S^2$ or $\RP^2$. For the case $M=S^2$ or $\RP^2$, a set of generators for the Brunnian braids on $M$ is given by our generating set together with the homotopy groups of a $2$-sphere.
\end{abstract}

\maketitle
\tableofcontents

\section{Introduction}
Let $M$ be a general compact connected surface, possibly with boundary components and
$B_n(M)$ denotes the $n$-strand braid group on a surface $M$. From the point
of view of braids compactness of a surface is not essential: braids rest the same if you replace a boundary component by a puncture. Essential is that the
number of punctures is finite, so that the fundamental group and the braid groups will be finitely generated.

A \textit{Brunnian braid} means a braid that becomes trivial after
removing any one of its strands. The formal definition of Brunnian
braids is given in Section~\ref{section2}. A typical example of $3$-strand Brunnian braid on a disk is the braid given by the expression $(\sigma_1^{-1}\sigma_2)^3$, where $\sigma_1$ and $\sigma_2$ are the standard
braid generators. In picture, $\sigma_1^{-1}\sigma_2$ corresponds to the braid operation on three strands which consists of taking the crossing of the first two strands with the first strand above the second one followed by the crossing of the last two strands with the last strand above the second one. This  is a usual operation that women are doing to construct their braids. 
If a woman repeats the braid operation $\sigma_1^{-1}\sigma_2$ $k$ times,
where $k$ is a multiple of $3$, then she obtains a Brunnian braid.

Let $\Brun_n(M)$ denote the set of the $n$-strand Brunnian braids. Then $\Brun_n(M)$ forms a subgroup of $B_n(M)$. A classical question proposed by G.~S.~Makanin~\cite{Makanin} in 1980 is to determine a set of
generators for Brunnian braids over the disk. Brunnian braids were called
\textit{smooth braids} by Makanin. This question was
answered by D.~L.~Johnson~\cite{Johnson} and G.~G.~Gurzo~\cite{Gurzo}.
A different approach to
this question can be found in~\cite{LW, Wu1}.
In 70ies H.~W.~Levinson \cite{Lev1, Lev2} defined a notion of {\it
$k$-decomposable} braid. It means a braid which becomes a trivial
after removal of any arbitrary $k$ strings. In his terminology
{\it decomposable} braid mean {\it $1$-decomposable} and so, Brunnian.

A connection between Brunnian braids and the homotopy groups of spheres
was given in~\cite{BCWW}. In particular, the following exact sequence
\begin{equation}
1\to \Brun_{n+1}(S^2)\to \Brun_n(D^2) \to \Brun_n(S^2)\to \pi_{n-1}(S^2)
\to 1
\label{eq:S2D2}
\end{equation}
was proved for $n>4$.

 In Birman's
book~\cite[Question 23, p. 219]{Birman}, she asked to determine a free basis for $\Brun_n(D^2)\cap R_{n-1}$ where
$$R_{n-1}=\Ker(B_n(D^2)\to B_n(S^2)).$$
Her motivation was that the kernel of the Gassner representation
is a subgroup of $\Brun_n(D^2)\cap R_{n-1}$.
From the exact sequence (\ref{eq:S2D2}) it follows that Birman's
question for $n>5$ is about a free basis of
Brunnian braids over the sphere $S^2$.
As far as we know this  question remains open.

The purpose of this article is to determine a set of generators for $\Brun_n(M)$ for a general surface $M$. We are able to determine a generating set for $\Brun_n(M)$ except two special cases where $M=S^2$ or $\RP^2$. For the case $M=S^2$ or $\RP^2$, we are able to determine a generating set for a (normal) subgroup of $\Brun_n(M)$ with the factor group given by $\pi_{n-1}(S^2)$.

We will use the notion of symmetric commutator subgroup. Given a
group $G$, and a set of its normal subgroups $R_1,\dots, R_n,\
(n\geq 2)$ denote
$$
[R_1,\dots, R_n]_S:=\prod_{\sigma\in
\Sigma_n}[[R_{\sigma(1)},R_{\sigma(2)}],\dots,R_{\sigma(n)}],
$$
where $\Sigma_n$ is the symmetric group of degree $n$.

Let $P_n(M)$ be the $n$-strand pure braid group on $M$. Let $D^2$ be a small disk in $M$. Then the inclusion $f\colon D\hookrightarrow M$ induces a group homomorphism
$$
f_*\colon P_n(D^2)\longrightarrow P_n(M).
$$
Recall that the pure Artin braid group $P_n(D^2)$ is generated by
the elements
$$
A_{i,j}=\sigma_{j-1}\sigma_{j-2}\cdots\sigma_{i+1}\sigma_i^2\sigma_{i+1}^{-1}\cdots\sigma_{j-2}^{-1}\sigma_{j-1}^{-1}
$$
for $1\leq i<j\leq n$. Let $A_{i,j}[M]=f_*(A_{i,j})$ and let $\la\la A_{i,j}[M]\ra\ra^P$ be the normal closure of $A_{i,j}[M]$ in $P_n(M)$. Note that a set of generators for $\la\la A_{i,j}[M]\ra\ra^P$ is given by $\beta A_{i,j}[M]\beta^{-1}$ for $\beta\in P_n(M)$. Thus a set of generators for the iterated subgroup
$$
[\la\la A_{1,n}[M]\ra\ra^P,\la\la A_{2,n}[M]\ra\ra^P,\ldots,
\la\la A_{n-1,n}[M]\ra\ra^P]_S
$$
can be given.

Now our determination for $\Brun_n(M)$ is as follows. For the cases $n\leq 3$, the determination is given in Propositions~\ref{proposition3.2}, ~\ref{proposition3.6} and ~\ref{brun3-rp2} by explicit computations. For $n\geq4$, the answer is as follows.

\begin{thm}[Theorem~\ref{theorem3.12}]\label{theorem1.1}
Let $M$ be a connected $2$-manifold and let $n\geq 4$. Let
$$
R_n(M)=[\la\la A_{1,n}[M]\ra\ra^P,\la\la
A_{2,n}[M]\ra\ra^P,\ldots, \la\la A_{n-1,n}[M]\ra\ra^P]_S
$$
be the symmetric commutator subgroup.
\begin{enumerate}
\item If $M\not=S^2$ or $\RP^2$, then
$$
\Brun_n(M)=R_n(M).
$$
\item If $M=S^2$ and $n\geq 5$, then there is a short exact sequence
$$
R_n(S^2)\hookrightarrow \Brun_n(S^2)\twoheadrightarrow \pi_{n-1}(S^2).
$$
\item If $M=\RP^2$, then there is a short exact sequence
$$
R_n(\RP^2)\hookrightarrow \Brun_n(\RP^2)\twoheadrightarrow \pi_{n-1}(S^2).
$$\end{enumerate}
\end{thm}

\noindent\textbf{Remarks.} \textbf{1.} Assertion 1 holds for each $n\geq 2$. Assertion 2 fails for $n=3, 4$. A free basis for $\Brun_4(S^2)$ has been given in~\cite{BCWW}. Assertion 3 fails for $n=2, 3$. The determination of $\Brun_2(\RP^2)$ and $\Brun_3(\RP^2)$ are given in section~\ref{section4}.

\textbf{2.} In the classical case where $M=D^2$, assertion 1 gives a better format for answering Makanin's question as we describe Brunnian braids as an explicit iterated commutator subgroup. Assertion 2 was essentially given in~\cite[Theorem 1.2]{BCWW}. Here we give an explicit determination for the kernel of $\Brun_n(S^2)\to \pi_{n-1}(S^2)$ for $n\geq5$. Assertion~3 gives a new connection between the Brunnian braids and the homotopy groups. The first case in assertion 3 ($n=4\sigma_2$)
is that the Hopf map $S^3\to S^2$ lifts to a $4$-strand Brunnian braid on $\RP^2$.

\textbf{3.} For the classical case the inclusion
$$
R_n(D^2)\hookrightarrow \Brun_n(D^2)
$$
was remarked by Levinson \cite{Lev2}, p.~53.

By Corollary~\ref{corollary2.4}, $\Brun_n(M)$ is a normal subgroup of $B_n(M)$ for $n\geq 3$. As an abstract  group, $\Brun_n(M)$ is a free group of infinite rank for $n\geq 3$ with $M\not=S^2$ or $\RP^2$, for $n\geq 5$ with $M=S^2$ and
for $n\geq 4$ with $M=\RP^2$. It comes out a natural question whether the factor group $B_n(M)/\Brun_n(M)$ is finitely presented. Our answer to this question is positive.

\begin{thm}\label{theorem1.2}
Let $M$ be a connected compact $2$-manifold. Then the factor groups $P_n(M)/\Brun_n(M)$ and $B_n(M)/\Brun_n(M)$ are finitely presented for each $n\geq 3$.
\end{thm}

The article is organized as follows. In Section~\ref{section2}, we give a review on Brunnian braids. The determination of a generating set for Brunnian braids is given in section~\ref{section3}, where Theorem~\ref{theorem3.12} is Theorem~\ref{theorem1.1}. In section~\ref{section4}, we compute the $3$-strand Brunnian braids on the projective plane. The proof of Theorem~\ref{theorem1.2} is given in section~\ref{section5}. In section~\ref{section6}, we give some remarks on Brunnian braids.

\section{Brunnian Braids}\label{section2}
\subsection{Configuration spaces and the braid groups}\label{subsection2.1}
Let $M$ be a topological space and let $M^n$ be the $n$-fold Cartesian product
of $M$. The \textit{$n$-th ordered configuration space} $F(M,n)$ is defined
by
$$
F(M,n)=\{(x_1,\ldots, x_n)\in M^n \ | \ x_i\not=x_j \textrm{ for }
i\not=j\}
$$
with subspace topology of $M^n$. The symmetric group $\Sigma_n$ acts
on $F(M,n)$ by permuting coordinates. The orbit space $$B(M,n)=F(M,n)/\Sigma_n$$ is called the \textit{$n$-th unordered configuration space}.
The \textit{braid group} $B_n(M)$ is defined to be the fundamental group
$\pi_1(F(M,n)/\Sigma_n)$. The \textit{pure braid group} $P_n(M)$ is defined to be the fundamental
group  of this pace $\pi_1(F(M,n)$. From the covering $F(M,n)\to F(M,n)/\Sigma_n$,
there is a short exact sequence of groups
$$
\{1\}\to P_n(M)\to B_n(M)\to \Sigma_n\to \{1\}.
$$

A geometric description of the elements in $B_n(M)$ can be given as follows.
Let $(q_1,\ldots,q_n)$ be the basepoint of $F(M,n)$ and let $$p\colon
F(M,n)\to F(M,n)/\Sigma_n$$ be the quotient map. The basepoint of
$F(M,n)/\Sigma_n$ is chosen to be $p(q_1,\ldots,q_n)$. Let
$[\lambda]$ be an element in $\pi_1(F(M,n)/\Sigma_n)$ represented by
a loop $\lambda\colon S^1\to F(M,n)/\Sigma_n$. Since $$p\colon
F(M,n)\to F(M,n)/\Sigma_n$$ is a covering, the loop $\lambda$ lifts
to a unique path $\tilde\lambda\colon [0,1]\to F(M,n)$ starting from
$\tilde \lambda(0)=(q_1,\ldots,q_n)$ and ending with $\tilde
\lambda(1)=(q_{\sigma(1)},\ldots, q_{\sigma(n)})$ for some
$\sigma\in\Sigma_n$. Let
$$
\tilde \lambda(t)=(\tilde\lambda_1(t),\ldots, \tilde\lambda_n(t))\in
F(M,n)\subseteq M^n.
$$
Then $\tilde \lambda_i(t)\not=\tilde\lambda_j(t)$ for $i\not=j$ and
any $0\leq t\leq 1$. The strands
$$\{(\tilde\lambda_i(t), t)\ | \ 1\leq i\leq n\}$$
in the cylinder $M\times [0,1]$ give the intuitive
braided description of $\lambda$. The precise definition of geometric braids are as follows.

Let $\{p_1,p_2,\ldots,p_n\}$ be $n$ distinct points in $M$. Consider the cylinder $M\times I$. A \textit{geometric braid}
$$
\rho=\{\rho_1,\ldots,\rho_n\}
$$
at the \textit{basepoints} $\{p_1,\ldots,p_n\}$ is a collection of $n$ paths in the cylinder $M\times I$ such that $\rho_i(t)= (\lambda_i(t), t)$ and
\begin{enumerate}
\item[1)] $\lambda_1(0)=p_1,\ldots, \lambda_n(0)=p_n$;
\item[2)] $\lambda_1(1)=p_{\sigma(1)},\ldots,\lambda_n(1)=p_{\sigma(n)}$ for some $\sigma\in \Sigma_n$;
\item[3)] $\lambda_i(t)\not=\lambda_j(t)$ for $0\leq t\leq 1$ and $i\not= j$;
\item[4)] each path $\lambda_i$ run monotonically with $t\in [0,1]$.
\end{enumerate}
Let $\rho=\{\rho_1,\ldots,\rho_n\}$ and $\rho'=\{\rho'_1,\ldots,\rho'_n\}$ be two geometric braids. We call $\rho$ is equivalent to $\rho'$, denoted by $\rho\sim \rho'$, if there exists a continuous sequence of geometric braids
$$\rho^s = (\lambda^s, t)=\{(\lambda_1^s(t),t),\
\ldots,(\lambda_n^s(t), t)\},\ 0\leq s\leq 1,$$ such that
\begin{enumerate}
\item[1)] $\lambda_1^s(0)=p_1,\ldots,\lambda_n^s(0))=p_n$ for each $0\leq s\leq 1$;
\item[2)] $\lambda_1^s(1)=
\lambda_1^0(1),\ldots,\lambda_n^s(1)=\lambda_n^0(1)$ for each $0\leq s\leq 1$;
\item[3)] $\lambda^0=\lambda$ and $\lambda^1=\lambda'$.
\end{enumerate}
In other words $\rho\sim \rho'$ if and only if they represent the same path homotopy class in the configuration space $F(M,n)$. A \textit{(geometric) braid} $\beta$ refers to a representative of the equivalence class of geometric braids.

The product of two geometric braids $\beta$ and $\beta'$ is defined to be the composition of the strands. More precisely, let $\beta$ be represented by $\rho=\{\rho_1,\ldots,\rho_n\}$ with $\rho_1(1)=p_{\sigma(1)},\ldots,\rho_n(1)=p_{\sigma(n)}$ and let $\beta'$ be represented by $\rho'=\{\rho'_1,\ldots,\rho'_n)$. Then the product $\beta\beta'$ is represented by
$$
\rho\ast\rho'=\{\rho_1\ast \rho'_{\sigma(1)},\ldots, \rho_n\ast \rho'_{\sigma(n)}\},
$$
where $\rho_i\ast\rho'_{\sigma(i)}$ is the path product.

\subsection{Removing Strands}\label{subsection2.2}
A \textit{simple (half-open) curve} in a space $M$ means a continuous injection $\theta\colon \mathbb{R}^+=[0,\infty)\to M$.
The distinct points $\{p_1,\ldots,p_n\}$ in $M$ are said \textit{well-ordered} with respect to a simple curve $\theta$ if there exists a sequence of points
on the interval
 $0\leq t_1<t_2<\cdots<t_n$ such that $p_i=\theta(t_i)$ for $1\leq i\leq n$.

Let $\mathbf{p}=\{p_1,\ldots,p_n\}$ and $\mathbf{p}'=\{p'_1,\ldots,p'_n\}$ be two sets of $n$ distinct well-ordered points with respect to $\theta$ with $p_i=\theta(t_i)$ and $p'_i=\theta(t'_i)$. Define
$$
L(\mathbf{p},\mathbf{p}')(s)=\{L(\mathbf{p},\mathbf{p}')_i(s)=\theta((1-s)t_i+st'_i) \ | \ 1\leq i\leq n\}
$$
for $0\leq s\leq 1$; $ L(\mathbf{p},\mathbf{p}')(s)\in M^n$. Observe that, for each $1\leq i<j\leq n$ and $0\leq s\leq 1$, $$(1-s)t_i+st'_i< (1-s)t_j+st'_j$$ as $t_i<t_j$ and $t'_i<t'_j$ and so $L(\mathbf{p},\mathbf{p}')(s)$ is a set of $n$ distinct well-ordered points with respect to $\theta$ for $0\leq s\leq 1$.

Now let  $\mathbf{p}=\{p_1,\ldots,p_n\}$ and $\mathbf{p}'=\{p'_1,\ldots,p'_n\}$ be two sets of $n$ distinct points on the curve $\theta$. There exist unique permutations $\sigma,\tau\in \Sigma_n$ such that
$$
\mathbf{p}_{\sigma}=\{p_{\sigma(1)},\ldots,p_{\sigma(n)}\}\textrm{ and } \mathbf{p}'_{\tau}=\{p'_{\tau(1)},\ldots,p'_{\tau(n)}\}
$$
are well-ordered with respect to $\theta$. We call
$$
L(\mathbf{p}_\sigma,\mathbf{p}'_{\tau})^{\sigma^{-1}}
=\{L(\mathbf{p}_\sigma,\mathbf{p}'_{\tau})_{\sigma^{-1}(i)} \ | \ 1\leq i\leq n\}
$$
an \textit{$n$-strand $\theta$-linear braid} from $\mathbf{p}$ to a permutation of $\mathbf{p}'$.


Let $M$ be a space with a simple curve $\theta$ and let the basepoints $\{p_1,p_2,\ldots, p_n \}$ of the braids on $M$ be well-ordered with respect to $\theta$. The \textit{system of removing strands} $d_i\colon B_n(M)\to B_{n-1}(M)$ is defined as follows:
\begin{enumerate}
\item[] Let $\beta\in B_n(M)$ be a braid represented by $\lambda=\{\lambda_1,\ldots,\lambda_n\}$ with
$$
\lambda_1(1)=p_{\sigma(1)},\ldots,\lambda_n(1)=p_{\sigma(n)}.
$$
Then the braid $d_i(\beta)$ is defined to be the equivalence class represented by the path product of the strands given by
$$
L\ast \{\lambda_1,\ldots, \lambda_{i-1},\lambda_{i+1},\ldots,\lambda_n\}\ast L',
$$
where $L$ is the $n-1$-strand $\theta$-linear braid from $\{p_1,\ldots,p_{n-1}\}$ to $\{p_1,\ldots, p_{i-1},p_{i+1},\ldots,p_n\}$ and $L'$ is the $n-1$-strand $\theta$-linear braid from $\{p_{\sigma(1)},\ldots, p_{\sigma(i-1)}, p_{\sigma(i+1)},\ldots,p_{\sigma(n)}\}$ to a permutation of $\{p_1,\ldots,p_{n-1}\}$.
\end{enumerate}
It follows from this definition that the operation $d_i$ does not
depend on the choice of $\lambda$ in the class $\beta$.
Intuitively, the operation $d_i\colon B_n(M)\to B_{n-1}(M)$ is obtained by forgetting the $i$-th strand and gluing back to the fixed choice of the basepoints using $\theta$-linear braids.

From now on we always assume that the space $M$ has a simple curve $\theta$ and the basepoints of the braids on $M$ are located on the curve $\theta$ starting with a set $\mathbf{p}$ of well-ordered points with respect to $\theta$ and ending with a permutation on $\mathbf{p}$. Recall that there is a short exact sequence
$$
1\to P_n(M)\to B_n(M)\to \Sigma_n\to 1.
$$
The braid group $B_n(M)$ acts by right on the letters $\{1,2,\ldots,n\}$ through the epimorphism $B_n(M)\to \Sigma_n$, which can be described as follows. Let $\beta$ be represented by an $n$-strand geometric braid
$$
\lambda=\{\lambda_i(t) \ | \ 1\leq i\leq n\}
$$
with $\lambda_i(0)=p_i$. Then $i\cdot \beta$ is given by the formula
$$
\lambda_i(1)=p_{i\cdot\beta}
$$
for $1\leq i\leq n$.

\begin{prop}~\cite[Proposition 4.2.1 (1)]{BCWW}\label{proposition2.1}
Let $M$ be a space with a simple curve. Then the operations $$d_i\colon B_n(M)\to B_{n-1}(M),\ 1\leq i\leq n,$$ satisfy the following identities:
\begin{enumerate}
\item[1)] $d_id_j=d_jd_{i+1}$ for $i\geq j$;
\item[2)] $d_i(\beta\beta')=d_i(\beta)d_{i\cdot \beta}(\beta')$.\hfill $\Box$
\end{enumerate}
\end{prop}

\smallskip

\noindent\textbf{Note.} In~\cite{BCWW}, the removing-strand operations are labeled by $d_0,\ldots, d_{n-1}$ for coinciding with simplicial terminology. The above identities are directly translated from~\cite[Proposition 4.2.1 (1)]{BCWW}.

\subsection{Brunnian Braids}\label{subsection2.3}

\begin{dfn}
Let $M$ be a space with a simple curve. A braid $\beta\in B_n(M)$ is called \textit{Brunnian} if $d_i(\beta)=1$ for each $1\leq i\leq n$. The set of $n$-strand Brunnian braids is denoted by $\Brun_n(M)$. For convention, any $1$-strand braid is regarded as a Brunnian braid.
\end{dfn}

Intuitively a Brunnian braid means a braid that becomes trivial after removing any one of its strands. If $\beta,\beta'\in\Brun_n(M)$, then
$$
d_i(\beta\beta')=d_i(\beta)d_{i\cdot\beta}(\beta')=1
$$
for $1\leq i\leq n$ and so the product $\beta\beta'\in \Brun_n(M)$. Similar $\beta^{-1}$ is Brunnian provided $\beta$ is. Thus $\Brun_n(M)$ is a subgroup of $B_n(M)$.

\begin{prop}\label{proposition2.3}
Let $M$ be a space with a simple curve. Then the subgroup $\Brun_n(M)\cap P_n(M)$ is normal in $B_n(M)$ for each $n\geq 1$.
\end{prop}
\begin{proof}
Let $\beta\in \Brun_n(M)\cap P_n(M)$ and let $\gamma\in B_n(M)$. Then
$$
\begin{array}{rcl}
d_i(\gamma\beta\gamma^{-1})&=&d_i(\gamma\beta)d_{i\cdot(\gamma\beta)}(\gamma^{-1})\\
&=&d_i(\gamma)d_{i\cdot\gamma}(\beta)d_{i\cdot(\gamma\beta)}(\gamma^{-1})\\
&=&d_i(\gamma)d_{i\cdot(\gamma\beta)}(\gamma^{-1})\\
\end{array}
$$
for $1\leq i\leq n$.
Since $\beta\in P_n(M)$, the elements $\gamma$ and $\gamma\beta$ have the same image in $\Sigma_n=B_n(M)/P_n(M)$ and so $i\cdot(\gamma\beta)=i\cdot\gamma$. The assertion follows from the equation
$$
1=d_i(1)=d_i(\gamma\gamma^{-1})=d_i(\gamma)d_{i\cdot\gamma}(\gamma^{-1})=d_i(\gamma)d_{i\cdot(\gamma\beta)}(\gamma^{-1}).
$$
\end{proof}

\begin{cor}\label{corollary2.4}
Let $M$ be a space with a simple curve. Then $\Brun_n(M)$ is a normal subgroup of $B_n(M)$ for $n\geq 3$.
\end{cor}
\begin{proof}
According to~\cite[Proposition 4.2.2]{BCWW}, $\Brun_n(M)\leq P_n(M)$ for $n\geq 3$ and hence the result.
\end{proof}

The case $n=2$ is exceptional for having $\Brun_n(M)$ to be normal in $B_n(M)$.

\begin{prop}
Let $M$ be a connected $2$-manifold. Then $\Brun_2(M)$ is a normal subgroup of $B_2(M)$ if and only if $\pi_1(M)=\{1\}$.
\end{prop}
\begin{proof}
If $\pi_1(M)=\{1\}$, then $B_2(M)=\Brun_2(M)$ as $B_1(M)=\pi_1(M)$.

Suppose that $\pi_1(M)\not=\{1\}$. Let $D^2$ be a small disk in $M\smallsetminus \partial M$. The inclusion $f\colon D^2\to M$ induces canonical maps $$(f,f)\colon F(D^2,2)\rightarrowtail F(M,2)\textrm{ and } (f,f)\colon F(D^2,2)/\Sigma_2\rightarrowtail F(M,2)/\Sigma_2.$$ Thus there is a commutative diagram of short exact sequences of groups
\begin{diagram}
1&\rTo &P_2(M)&\rTo & B_2(M)&\rTo &\Sigma_2&\rTo &1\\
 &  &\uTo& &\uTo>{(f,f)_*}&&\uEq&&\\
1&\rTo &P_2(D^2)&\rTo & B_2(D^2)&\rTo &\Sigma_2&\rTo &1.\\
\end{diagram}
Let $\sigma_1$ be a generator for $B_2(D^2)=\Z$. Then $(f,f)_*(\sigma_1)\not=1$ in $B_2(M)$ as it has the nontrivial image in $\Sigma_2=B_2(M)/P_2(M)$. From the commutative diagram
\begin{diagram}
B_2(D^2)&\rTo^{(f,f)_*}&B_2(M)\\
\dTo>{d_i}&&\dTo>{d_i}\\
B_1(D^2)=\{1\}&\rTo^{f_*}&B_1(M)\\
\end{diagram}
for $i=1,2$, the element $\beta=(f,f)_*(\sigma_1)$ is a Brunnian braid on $M$. Let $p_1$ be the basepoint of $M$. Choose a loop
$$
\omega\colon [0,1]\to M
$$
with $\omega(0)=\omega(1)=p_1$ representing a nontrivial element in $\pi_1(M)$. Take the second basepoint $p_2$ such that $p_2$ is not on the curve $\omega([0,1])$ and construct a $2$-strand braid $\gamma$ represented by
$$
\rho(t)=\{\rho_1(t),\rho_2(t)\}
$$
with $\rho_1(t)=(\omega(t), t)$ and $\rho_2(t)=(p_2, t)$ for $0\leq t\leq 1$ in the cylinder $M\times I$. Then $d_1(\gamma)=\{1\}$ as represented by the straight line-segment given by $\rho_2$, and $d_2(\gamma)=[\omega]\not=1$ the path homotopy class represented by $\omega$. Observe that $\gamma$ is a pure braid. We have $d_i(\gamma^{-1})=(d_i(\gamma))^{-1}$. From
$$
\begin{array}{rcl}
d_1(\gamma\beta\gamma^{-1})&=&d_1(\gamma)d_{1\cdot \gamma}(\beta)d_{1\cdot(\gamma\beta)}(\gamma^{-1})\\
&=&d_1(\gamma)d_{1}(\beta)d_{2}(\gamma^{-1})\\
&=&d_1(\gamma)d_{1}(\beta)d_{2}(\gamma)^{-1}\\
&=&1\cdot 1\cdot [\omega]^{-1}\\
&\not=&1,\\
\end{array}
$$
the conjugation $\gamma\beta\gamma^{-1}$ is not Brunnian and so $\Brun_2(M)$ is not normal. This finishes the proof.
\end{proof}

\section{Generating Sets for Brunnian Braids on Surfaces}\label{section3}
In this section, $M$ is a connected compact $2$-dimensional (oriented or non-oriented) manifold. The classical Fadell-Neuwirth Theorem will be useful in computations.

\begin{thm}\cite{FN}\label{theorem3.1}
The coordinate projection
$$
\delta^{(i)}\colon F(M,n)\to F(M,n-1), \ (x_1,\ldots,x_n)\mapsto (x_1,\ldots,x_{i-1},x_{i+1},\ldots,x_n)
$$
is a fiber bundle with fiber $M\smallsetminus Q_{n-1}$, where $Q_{n-1}$ is a set of $(n-1)$ distinct points in $M$.\hfill $\Box$
\end{thm}
\begin{prop} Up to the change of base-point for the pure braid group
$P_n(M)$ the homomorphism $d_i$ coincides with homomorphism of
fundamental groups induced by $\delta^{(i)}$:
$$
d_i =\delta^{(i)}_* h_i \colon P_n(M)\to P_{n-1}(M)
$$
where $h_i$ is the automorphism of $\pi_1 (F(M,n-1))$ induced by the
change of base-points
$$(F(M,n-1), (p_1, \dots, p_{i-1}, p_{i+1}, \dots, p_n)) \to
 (F(M,n-1), (p_1, \dots, p_{n-1})).
 $$
 \hfill $\Box$
\end{prop}

Let $D^2$ be a small disk in $M\smallsetminus \partial M$. The basepoints $\{p_1,p_2,\ldots\}$ for the braids on $M$ are chosen inside $D^2\smallsetminus\partial D^2$. The embedding $f\colon D^2\rightarrowtail M$ induces a map
$$
f^n\colon F(D^2,n)/\Sigma_n \rightarrowtail F(M,n)/\Sigma_n
$$
and so a group homomorphism
$$
f^n_*\colon B_n(D^2)=\pi_1(F(D^2,n)/\Sigma_n) \longrightarrow B_n(M)=\pi_1(F(M,n)/\Sigma_n)
$$
with a commutative diagram
\begin{diagram}
B_n(D^2)&\rTo^{f^n_*}&B_n(M)\\
\dOnto&&\dOnto\\
B_n(D^2)/P_n(D^2)=\Sigma_n&\rEq&\Sigma_n=B_n(M)/P_n(M).\\
\end{diagram}
For any braid $\beta\in B_n(D^2)$, we write $\beta[M]$ (or simply $\beta$ if there are no confusions) for the braid $f^n_*(\beta)$ on $M$.

Recall that the Artin braid group $B_n(D^2)$ is generated by $\sigma_1,\ldots,\sigma_{n-1}$ with defining relations
\begin{enumerate}
\item $\sigma_i\sigma_j=\sigma_j\sigma_i$ for $|i-j|\geq 2$ and
\item $\sigma_i\sigma_{i+1}\sigma_i=\sigma_{i+1}\sigma_i\sigma_{i+1}$ for each $i$,
\end{enumerate}
where as a geometric braid, $\sigma_i$ is the canonical $i$th elementary braid of $n$-strands that twists the positions $i$ and $i+1$ once and puts the trivial strands on the remaining positions. Also recall that the pure Artin braid group $P_n(D^2)$ is generated by
$$
A_{i,j}=\sigma_{j-1}\sigma_{j-2}\cdots\sigma_{i+1}\sigma_i^2\sigma_{i+1}^{-1}\cdots\sigma_{j-2}^{-1}\sigma_{j-1}^{-1}
$$
for $1\leq i<j\leq n$.

\subsection{$2$-strand Brunnian Braids}\label{subsection3.1}
\begin{prop}\label{proposition3.2}
Let $M$ be any connected $2$-manifold. Then the $2$-strand Brunnian braids are determined as follows
\begin{enumerate}
\item[1)] $\Brun_2(M)\cap P_2(M)$ is the normal closure of the element $A_{1,2}$ in $B_2(M)$.
\item[2)] $\Brun_2(M)$ is the subgroup of $B_2(M)$ generated by $\Brun_2(M)\cap P_2(M)$ and $\sigma_1$: $\Brun_2(M) =\la \Brun_2(M)\cap P_2(M),\sigma_1\ra$.
\end{enumerate}
\end{prop}
\begin{proof}
(1) Let $\la\la A_{1,2}\ra\ra$ be the normal closure of $A_{1,2}$ in $B_2(M)$. By Proposition~\ref{proposition2.3}, $\Brun_2(M)\cap P_2(M)$ is normal in $B_2(M)$. Since $A_{1,2}$ is a pure Brunnian braid,
$$
\la\la A_{1,2}\ra\ra\leq \Brun_2(M)\cap P_2(M).
$$
To see the equality, consider the commutative diagram of fiber sequences
\begin{diagram}
F&\rTo& M\smallsetminus\{p_2\}&\rInto^{i'} &M\\
\dTo&&\dInto>{i_1}&&\dEq\\
M\smallsetminus\{p_1\}&\rInto^{i_2}&F(M,2)&\rTo^{d_2}&M\\
\dInto>{i}&&\dTo>{d_1}&&\dTo\\
M&\rEq&M&\rTo&\ast,\\
\end{diagram}
where $i_2(x)=(p_1,x)$ and $i_1(x)=(x,p_2)$ and $F$ is a homotopy fiber
of $i$, what is equivalent to a fiber of $i'$.
From the middle row, there is an exact sequence
\begin{multline}
\pi_2(M)\rTo \pi_1(M\smallsetminus \{p_1\})\rTo^{i_{2*}} \pi_1(F(M,2))=\\
P_2(M)\rOnto^{d_2}
 \pi_1(M)=P_1(M).
 \label{ex_se_pi}
\end{multline}
Note that
\begin{multline*}
\Brun_2(M)\cap P_2(M)=\\
\Ker(d_1\colon P_2(M)\to P_1(M))\cap \Ker(d_2\colon P_2(M)\to P_1(M)).
\end{multline*}
Consider the following  diagram of the short exact sequences of groups
\begin{equation}
\begin{diagram}
\la\la \omega\ra\ra &\rOnto^{i_{2*}|}&\Brun_2(M)\cap P_2(M)\\
\dInto&&\dInto\\
\pi_1(M\smallsetminus\{p_1\})&\rOnto^{i_{2\ast}}&\Ker(d_2\colon P_2(M)\to P_1(M))\\
\dOnto>{i_*}&&\dTo>{d_1}\\
P_1(M)&\rEq&P_1(M),\\
\end{diagram} \label{omega}
\end{equation}
where $\omega\in \pi_1(M\smallsetminus\{p_1\})$ is represented by a small circle around $p_1$. Its commutativity follows from construction and
epimorphisms follow from the exact sequence (\ref{ex_se_pi}).
It follow from the diagram (\ref{omega}) that
 $\Brun_2(M)\cap P_2(M)$ is the normal closure of $i_{2*}(\omega)$ in $\Ker(d_2)$. From the commutative diagram,
\begin{diagram}
 \pi_1(M\smallsetminus \{p_1\})&\rTo^{i_{2*}} &\pi_1(F(M,2))&\rOnto^{d_2} &\pi_1(M)\\
 \uTo>{f_*}&&\uTo>{(f\times f)_*}&&\uTo\\
 \pi_1(D^2\smallsetminus \{p_1\})=\Z&\rTo^{i_{2*}}_{\cong}& \pi_1(F(D,2))&\rOnto^{d_2} &\pi_1(D^2)=1,\\
\end{diagram}
we get
$$
i_{2*}(\omega)=A_{1,2}^{\pm 1}
$$
and hence it follows assertion 1).

(2) Note  that $\sigma_1$ is Brunnian. From  the short exact sequence
$$
1\to P_2(M)\to B_2(M)\to \Sigma_2\to 1.
$$
we get the following commutative diagram
\begin{diagram}
1&\to \Brun_2(M)\cap P_2(M)&\to \Brun_2(M)\to &\Sigma_2&\to 1\\
&\dInto&\dInto&\dEq\\
1&\longrightarrow  P_2(M)&\longrightarrow \ \ B_2(M)\to &\Sigma_2&\to 1.\\
\end{diagram}
 and the assertion follows.
\end{proof}
\begin{cor}
Let $M$ be a connected $2$-manifold. Then $$B_2(M)/(\Brun_2(M)\cap P_2(M))$$ is the quotient group of $B_2(M)$ by adding the single relation $$A_{1,2}=\sigma_1^2=1.
 $$
 \hfill $\Box$
\end{cor}

\subsection{Homotopy properties of configuration spaces of surfaces}\label{subsection3.1,5}
The following (well-known) fact will be useful for the computations in next subsections.
\begin{lem}\label{lemma3.5}
Let $M$ be a connected $2$-manifold.
\begin{enumerate}
\item If $M\not=S^2$ or $\RP^2$, then $F(M,n)$ is a $K(\pi,1)$-space for $n\geq 1$. In particular, $\pi_2(F(M,n))=0$ for $n\geq 1$.
\item $\pi_2(F(S^2,n))=0$ for $n\geq 3$.
\item $\pi_2(F(\RP^2,n))=0$ for $n\geq 2$.
\end{enumerate}
\end{lem}
\begin{proof}
Assertion~(1) follows from the fact that $M$ and $M\setminus Q_{n-1}$
are $K(\pi,1)$ spaces together with Fadell-Neuwirth fibration (Theorem~\ref{theorem3.1}).

To prove the assertion~(2) we note at first that $F(S^2,3)$ is homotopy
equivalent to the rotation group $SO(3)$ by orthogonalization process:
$F(S^2,3)\simeq SO(3)$ and the rotation group is homeomorphic to
the real projective space $\RP^3$ (classical fact). So,
$\pi_2(F(S^2,3))=0$.
There is also the following   decomposition formula~\cite[Corollary 2.3]{BCWW}
$$
F(S^2,n)\simeq F(\mathbb{R}^2\smallsetminus\{0,1\},n-3)\times F(S^2,3)
$$
for $n\geq 4$, from which the assertion~(2) follows.

(3). Let us consider  the following exact sequence
\begin{multline}
\pi_2(\RP^2\smallsetminus\{p_1\})=0\to \pi_2(F(\RP^2,2))\to \pi_2(\RP^2)\rTo^{\partial} \\
\pi_1(\RP^2\smallsetminus\{p_1\})
\to
\pi_1(F(\RP^2,2))\to \pi_1(\RP^2)\to 1.
\label{extens}
\end{multline}
We note that the homomorphism $\partial$ is a monomorphism between
two groups isomorphic to $\Z$, actually it is equal to multiplication
by $4$ as $\pi_1(F(\RP^2,2))= P_2(\RP^2) = {\bold Q}_8$, the quaternion group
\cite{VB} and
$\pi_1(\RP^2)=\Z/2$. So,
$\pi_2(F(\RP^2,2))=0$. For $n\geq 3$, consider the fiber sequence
$$
F(\RP^2\smallsetminus Q_2,n-2)\rInto F(\RP^2,n)\rOnto F(\RP^2,2),
$$
The assertion follows from the facts that $F(\RP^2\smallsetminus Q_2,n-2)$ is a $K(\pi,1)$-space by assertion 1 and $\pi_2(F(\RP^2,2))=0$.
\end{proof}

\subsection{$3$-strand Brunnian Braids}\label{subsection3.2}
Now we are going to determine the $3$-strand Brunnian braids on $M$. By~\cite[Proposition 4.2.2]{BCWW}, $$\Brun_n(M)\subseteq P_n(M)$$ for $n\geq 3$. Thus the determination is given by
$$
\Brun_n(M)=\Brun_n(M)\cap P_n(M)=\bigcap_{i=1}^n\Ker(d_i\colon P_n(M)\to P_{n-1}(M))
$$
for $n\geq 3$.

For a subset $S$ in $P_n(M)$, we write $\la\la S\ra\ra^P$ for the normal closure of $S$ in $P_n(M)$ while we keep the notation $\la\la S\ra\ra$ for the normal closure of $S$ in $B_n(M)$.

\begin{prop}\label{proposition3.6}
Let $M$ be a connected $2$-manifold. Then the $3$-strand Brunnian braids on $M$ are determined as follows:
\begin{enumerate}
\item[1)] $\Brun_3(S^2)=P_3(S^2)=\Z/2$.
\item[2)] For $M\not=S^2$ or $\RP^2$,
$$
\Brun_3(M)=[\la\la A_{1,3}\ra\ra^P,\la\la A_{2,3}\ra\ra^P]
$$
the commutator subgroup of the normal closures in $P_3(M)$ generated by $A_{1,3}$ and $A_{2,3}$, respectively.
\end{enumerate}
\end{prop}
\begin{proof}
Assertion 1 follows directly from the fact that $P_3(S^2)=\Z/2$ and $P_2(S^2)=\{1\}$. For assertion~$2$, observe that $d_kA_{i,j}=1$ for $k=i,j$. Thus
$$
\la\la A_{i,3}\ra\ra^P\leq \Ker(d_3\colon P_3(M)\to P_2(M))\cap \Ker(d_i\colon P_3(M)\to P_2(M))
$$
for $i=1,2$ and so, the inclusion
$$
[\la\la A_{1,3}\ra\ra^P,\la\la A_{2,3}\ra\ra^P]\leq \Brun_3(M)
$$
is simple.


From the commutative diagram of the fiber sequences
\begin{equation}
\begin{diagram}
M\smallsetminus\{p_1,p_2\}&\rInto^{i_3} &F(M,3)&\rTo^{\delta_3}& F(M,2)\\
\dInto&&\dTo>{\delta_2}&&\dTo>{\delta_2}\\
M\smallsetminus\{p_1\}&\rInto^{i'_2}& F(M,2)&\rTo^{\delta_2}& M,\\
\end{diagram}
\label{equation3.1}
\end{equation}
where $i_3(x)=(p_1,p_2,x)$ and $i'_2(x)=(p_1,x)$, together with the facts that $\pi_2(M)=0$ and $\pi_2(F(M,2))=0$ (Lemma~\ref{lemma3.5}), there is a commutative diagram of short exact exact sequences
\begin{equation}\label{equation3.2}
\begin{diagram}
\pi_1(M\smallsetminus\{p_1,p_2\})&\rInto^{{i_3}_*} &P_3(M)&\rOnto^{d_3}& P_2(M)\\
\dTo>{d_2|}&&\dTo>{d_2}&&\dTo>{d_2}\\
\pi_1(M\smallsetminus\{p_1\})&\rInto^{{i'_2}_*}&P_2(M)&\rOnto^{d_2}& P_1(M).\\
\end{diagram}
\end{equation}
It follows from this diagram that
$$
{i_3}_*\colon \Ker(d_2|)\longrightarrow \Ker(d_3\colon P_3(M)\to P_2(M))\cap \Ker(d_2\colon P_3(M)\to P_2(M))
$$
is an isomorphism. Since $d_2|\colon \pi_1(M\smallsetminus\{p_1,p_2\})\to \pi_1(M\smallsetminus\{p_1\})$ is induced by the inclusion
$$
M\smallsetminus \{p_1,p_2\}\rInto M\smallsetminus\{p_1\},
$$
$\Ker(d_2|)$ is the normal closure in $\pi_1(M\smallsetminus\{p_1,p_2\})$ generated by $[\omega_2]$, where $\omega_2$ is a small circle around $p_2$.
Similarly, the inclusion $$M\smallsetminus\{p_1,p_2\}\rInto M\smallsetminus\{p_2\}$$ induces a homomorphism
$$
d_1|\colon \pi_1(M\smallsetminus \{p_1,p_2\})\rTo \pi_1(M\smallsetminus\{p_2\})
$$
with the property that
$$
{i_3}_*\colon \Ker(d_1|)\longrightarrow \Ker(d_3\colon P_3(M)\to P_2(M))\cap \Ker(d_1\colon P_3(M)\to P_2(M)).
$$
is an isomorphism and $\Ker(d_1|)$ is the normal closure in $\pi_1(M\smallsetminus\{p_1,p_2\})$ generated by the homotopy class
$[\omega_1]$, where $\omega_1$ is a small circle around $p_1$. Thus
\begin{equation}
{i_3}_*\colon \Ker(d_1|)\cap \Ker(d_2|)\longrightarrow \Brun_3(M)
\label{iso_i_3}
\end{equation}
is an isomorphism. By applying the Brown-Loday Theorem~\cite{Brown-Loday} to the homotopy push-out diagram of $K(\pi,1)$-spaces
\begin{diagram}
M\smallsetminus\{p_1,p_2\}&\rInto& M\smallsetminus\{p_1\}\\
\dInto&&\dInto\\
M\smallsetminus\{p_2\}&\rInto&M,\\
\end{diagram}
we get an isomorphism
$$
\frac{\Ker(d_1|)\cap \Ker(d_2|)}{[\Ker(d_1|),\Ker(d_2|)]}\cong \pi_2(M)=0
$$
and so
$$
\Ker(d_1|)\cap \Ker(d_2|)=[\Ker(d_1|),\Ker(d_2|)].
$$
Together with the isomorphism (\ref{iso_i_3}) this gives
\begin{equation}\label{equation3.3}
\Brun_3(M)=
[\la\la {i_3}_*([\omega_1])\ra\ra^P,\la\la {i_3}_*([\omega_2])\ra\ra^P].
\end{equation}
Note that the basepoints $\{p_1,p_2\}$ are chosen in the interior of the small disk $D^2$. From the commutative diagram
\begin{equation}\label{equation3.4}
\begin{diagram}
\pi_1(M\smallsetminus\{p_1,p_2\})&\rInto^{{i_3}_*} &P_3(M)&\rOnto^{d_3}& P_2(M)\\
\uTo>{f_*}&&\uTo>{{f}^3_*}&&\uTo>{{f}^2_*}\\
\pi_1(D^2\smallsetminus\{p_1,p_2\})&\rInto^{{i_3}_*}&P_3(D^2)&\rOnto^{d_3}& P_2(D^2),\\
\end{diagram}
\end{equation}
we have ${i_3}_*([\omega_1])=A_{1,3}^{\pm1}$ and ${i_3}_*([\omega_2])=A_{2,3}^{\pm1}$. Assertion $2$ follows from the replacing ${i_3}_*([\omega_i])$ by $A_{i,3}$ in equation~(\ref{equation3.3}).

\end{proof}

\subsection{Colimits of classifying spaces} Given a group $G$ and its normal subgroups
$R_1,\dots, R_n$, let us define their {\it complete commutator
subgroup} as follows
\begin{equation}
[[R_1,R_2,\ldots,R_n]]:=\prod\limits_{\begin{array}{c}
I\cup J=\{1,2,\ldots,n\}\\
I\cap J=\emptyset\\
\end{array}}\left[\bigcap\limits_{i\in I} R_i, \bigcap\limits_{j\in J}R_j\right].
\label{eq:comp_com}
\end{equation}
It is clear that $$ [[R_1,\dots, R_n]]\subseteq R_1\cap
\dots \cap R_n
$$
and that the quotient
$$
\frac{R_1\cap \dots\cap R_n}{[[R_1,\dots, R_n]]}
$$
is an abelian group with a natural structure of $\mathbb
Z[G/R_1\dots R_n]$-module, where the action is defined via
conjugation in $G$. An $n$-tuple of normal subgroups $(R_1,\dots,
R_n)$ is called {\it connected} in $G$ if either $n\leq 2,$ or
$n\geq 3$ and for all subsets $I, J \subset \{1,\cdots,n\}$ with $
|I|\ge 2, |J|\ge 1$ (without conditions of formula (\ref{eq:comp_com}))
the following equality holds:
\begin{equation}\label{connectivity}\left( \bigcap_{i\in I} R_i \right) \left( \prod_{j\in J}R_j\right) =
\bigcap_{i\in I} \left( R_i (\prod_{j\in J}R_j)
\right)\,.\end{equation}

We will make use of the following result from
\cite{Ellis-Mikhailov}

\begin{thm}\label{ellis-mikhailov}
Let $G$ be a group, $n\geq 2,$ and $(R_1,\dots, R_n)$ an $n$-tuple
of normal subgroups in $G$ such that the $(n-1)$-tuples
$(R_1,\dots, \hat R_i,\dots, R_n)$ are connected for all $1\leq
i\leq n$. Let $X$ be the topological space arising as the colimit
of classifying spaces $K(G/\prod_{i\in I}R_i,1)$, where $I$ ranges
over all subsets $I \subsetneq \{1,\dots,n\}$. Then there is an
isomorphism of abelian groups
$$\pi_n(X)\simeq \frac{R_1\cap \dots \cap R_n}{[[R_1,\dots,R_n]]}.
$$
\end{thm}

\subsection{$n$-strand Brunnian Braids for $n\geq 4$}\label{subsection3.3}
Now we are going to determine $\Brun_n(M)$ for $n\geq 4$. The case $\Brun_4(S^2)$ has been determined in~\cite[Proposition 7.2.2]{BCWW}. Our computation will exclude this special case.

\begin{lem}\label{lemma3.10}
Let $M$ be a connected $2$-manifold. Let
$$
d_k\colon P_n(M)\to P_{n-1}(M)
$$
be the operation that removes the $k$th strand.
\begin{enumerate}
\item Suppose that $M\not=S^2$ or $\RP^2$. Then, for $n\geq 2$,
$$
\Ker(d_n)\cap\Ker(d_k)=\la\la A_{k,n}\ra\ra^P
$$
for $1\leq k\leq n-1$ and therefore
$$
\Brun_n(M)=\bigcap_{k=1}^{n-1} \la\la A_{k,n}\ra\ra^P.
$$
Moreover ${i_n}_*\colon \pi_1(M\smallsetminus\{p_1,p_2,\ldots, p_{n-1}\})\to P_n(M)$ is a monomorphism with
$$
{i_n}_*(\Ker(d_k|))=\la\la A_{k,n}\ra\ra^P,$$
where $$
d_k|\colon \pi_1(M\smallsetminus \{p_1,\ldots,p_{n-1}\})\longrightarrow \pi_1(M\smallsetminus\{p_1,\ldots,p_{k-1},p_{k+1},\ldots,p_{n-1}\})
$$
is the group homomorphism induced by the inclusion by filling back the missing point $p_k$.
\item If $M=S^2$, then the above statement holds for $n\geq 5$.
\item If $M=\RP^2$, then the above statement holds for $n\geq 4$. \hfill $\Box$
\end{enumerate}
\end{lem}
\begin{proof}
 Diagram~(\ref{equation3.1}) can be extended to general case and so we have the starting commutative diagram
\begin{equation}\label{equation3.5}
\begin{diagram}
\pi_1(M\smallsetminus\{p_1,p_2,\ldots, p_{n-1}\})&\rInto^{{i_n}_*} &P_n(M)&\rOnto^{d_n}& P_{n-1}(M)\\
\dTo>{d_k|}&&\dTo>{d_k}&&\dTo>{d_{n-1}}\\
\pi_1(M\smallsetminus\{p_1,\ldots,p_{k-1},p_{k+1},\ldots,p_{n-1}\})&\rInto^{{i'_{n-2}}_*}&P_{n-1}(M)&\rOnto^{d_{n-2}}& P_{n-2}(M)\\
\end{diagram}
\end{equation}
for $1\leq k\leq n-1$. If $M\not=S^2$ or $n>4$, then
$$
{{i'_{n-2}}_*}\colon \pi_1(M\smallsetminus\{p_1,\ldots,p_{k-1},p_{k+1},\ldots,p_{n-1}\})\to P_{n-1}(M)
$$
is a monomorphism because $\pi_2(F(M,n-2))=0$ for $n\geq 4$ with $M\not=S^2$ or $n\geq 5$ with $M=S^2$ (Lemma~\ref{lemma3.5}).
Thus
$$
i_{n*}\colon \Ker(d_k|)\longrightarrow \Ker(d_n)\cap \Ker (d_k)
$$
is an isomorphism if $M\not=S^2$ or $n>4$. Note that $\Ker(d_k|)$ is the normal closure in $\pi_1(M\smallsetminus\{p_1,p_2,\ldots, p_{n-1}\})$ generated by the homotopy class
$[\omega_k]$, where $\omega_k$ is a small circle around $p_k$. By the same reasons as in diagram~(\ref{equation3.4}), we have $i_{n\ast}([\omega_k])=A_{k,n}^{\pm}$ and so
$$
\Ker(d_n)\cap \Ker (d_k)=i_{n\ast}(\Ker(d_k|))\leq \la\la A_{k,n}\ra\ra^P.$$
On the other hand, $\la\la A_{k,n}\ra\ra^P\leq \Ker(d_n)\cap \Ker (d_k)$ because $A_{k,n}$ lies in the normal subgroup $\Ker(d_n)\cap \Ker (d_k)$. Thus, in the case $M\not=S^2$ or $n>4$,
$$
\Ker(d_n)\cap \Ker (d_k)=\la\la A_{k,n}\ra\ra^P$$
and hence the result.
\end{proof}

The remaining question is of course how to determine the
intersection of the normal subgroups $\la\la A_{k,n}\ra\ra^P$ for
general $n$.

\begin{thm}\label{theorem3.11}
Let $M$ be a connected $2$-manifold and let $\{p_1,\ldots,p_n\}$ be the set of $n$ distinct points in $M\smallsetminus\partial M$. Let
$$
d_i|\colon \pi_1(M\smallsetminus \{p_1,\ldots,p_n\})\longrightarrow \pi_1(M\smallsetminus\{p_1,\ldots,p_{i-1},p_{i+1},\ldots,p_n\})
$$
be the group homomorphism induced by the inclusion by filling back the missing point $p_i$. Then
$$
\left(\bigcap_{i=1}^n\Ker(d_i|)\right)/[\Ker(d_1|),
\Ker(d_2|),\ldots,\Ker(d_n|)]_S\cong \pi_n(M)
$$
for each $n\geq 2$.
\end{thm}
\begin{proof}
Observe that the surface $M$ can be viewed as a colimit of the
spaces $M\smallsetminus \sqcup_{i\in I}p_i,$ where $I$ ranges over
all subsets $I\subsetneq \{1,\dots, n\}.$ Denote
$G:=\pi_1(M\smallsetminus \{p_1,\ldots,p_n\})$ and
$R_i:=\Ker(d_i|)$. Since punctured surfaces are aspherical, the
spaces $M\smallsetminus \sqcup_{i\in I}p_i$ are classifying spaces
for groups $G/\prod_{i\in I}R_i$. Let us check that the
connectivity condition (\ref{connectivity}) holds for every
$(n-1)$-tuple of subgroups $(R_1,\dots, \hat R_m,\dots, R_n),\
1\leq m\leq n$. For $n=2,3,$ the connectivity condition holds by
definition. We prove the statement by induction on $n$. We fix
number $m$: $1\leq m\leq n$ and prove the connectivity
(\ref{connectivity}) of the $(n-1)$-tuple $(R_1,\dots, \hat
R_m,\dots, R_n).$ Let $I, J\subseteq \{1,\dots, \hat m,\dots,
n\}.$
 Suppose
that $I\cap J\neq \emptyset$. Then the left and right hand sides
of (\ref{connectivity}) are equal to $\prod_{j\in J}R_j$ and the
condition is proved. So, we can assume that $I\cap J=\emptyset$.
Consider the epimorphism
$$
f_J: G\to G/\prod_{j\in J}R_j.
$$
The condition (\ref{connectivity}) is equivalent to the condition
\begin{equation}\label{conn1}
f_J(\bigcap_{i\in I}R_i)= \bigcap_{i\in I} f_J(R_i).
\end{equation} Any punctured surface has a free
fundamental group and
$$
f_J(R_i)=\Ker\{\pi_1(M\smallsetminus \sqcup_{k\in I}p_k)\to
\pi_1(M\smallsetminus \sqcup_{k\in I,\ k\neq i}p_k)\}.
$$
By induction we have $$\bigcap_{i\in I}R_i=[[R_{i_1},\dots,
R_{i_{|I|}}]]$$ for $I=\{i_1,\dots, i_{|I|}\}$ due to Theorem
\ref{ellis-mikhailov} and the fact that punctured surface is
aspherical. The same argument shows that
$$
\bigcap_{i\in I} f_J(R_i)=[[f_J(R_{i_1}),\dots, f_J(R_{i_{|I|}})]]
$$
(we repeat argument for the punctured surface with discs glued to
$|J|$ boundary components, the surface remains punctured since
$M\smallsetminus{p_1,\dots, p_n}$ has at least $n$ boundary
components). The same argument shows that
\begin{align*}
& [[R_{i_1},\dots, R_{i_{|I|}}]]=[R_{i_1},\dots, R_{i_{|I|}}]_S\\
& [[f_J(R_{i_1}),\dots, f_J(R_{i_{|I|}})]]=[f_J(R_{i_1}),\dots,
f_J(R_{i_{|I|}})]_S
\end{align*}
Since $f_J$ is a homomorphism, the condition (\ref{conn1}) and
hence (\ref{connectivity}) follow. Again observe that
$$
[[R_1,\dots, R_n]]=[R_1,\dots, R_n]_S,
$$
hence the needed statement follows from
Theorem~\ref{ellis-mikhailov}.
\end{proof}

\begin{thm}\label{theorem3.12}
Let $M$ be a connected $2$-manifold and let $n\geq 4$. Let
$$
R_n(M)=[\la\la A_{1,n}[M]\ra\ra^P,\la\la
A_{2,n}[M]\ra\ra^P,\ldots, \la\la A_{n-1,n}[M]\ra\ra^P]_S
$$
be the symmetric commutator subgroup.
\begin{enumerate}
\item If $M\not=S^2$ or $\RP^2$, then
$$
\Brun_n(M)=R_n(M).
$$
\item If $M=S^2$ and $n\geq 5$, then there is a short exact sequence
$$
R_n(S^2)\rInto \Brun_n(S^2)\rOnto \pi_{n-1}(S^2).
$$
\item If $M=\RP^2$, then there is a short exact sequence
$$
R_n(\RP^2)\rInto \Brun_n(\RP^2)\rOnto \pi_{n-1}(S^2).
$$\end{enumerate}
\end{thm}
\begin{proof}
By Lemma~\ref{lemma3.10},
$$
\Brun_n(M)=\bigcap_{i=1}^n\la\la A_{i,n}\ra\ra^P
$$
and $\la\la A_{k,n}\ra\ra^P=i_*(\Ker(d_k|))$. The assertion follows by Theorem~\ref{theorem3.11}
\end{proof}

\section{$3$-strand Brunnian Braids on the Projective Plane}\label{section4}

\subsection{Braid group of the projective plane}

There exist several presentations of the group
$B_n(\RP^2)$, see, for example \cite{VB, Gon-Gu}. We will use a presentation  similar to presentations of the surface braid group from
\cite{VB}.
\begin{thm}
 The group  $B_n(\RP^2)$ can be presented as
having  the set of generators
$$
 \sigma_1, \sigma_2, \ldots, \sigma_{n-1}, \rho,
$$
where in the braid $\rho$
the first string represents a nontrivial element of the fundamental
group and the rest of the braid is trivial; the generators
 $\sigma_1, \sigma_2, \ldots, \sigma_{n-1}$ are the images of classical braid generators
of the disk; the set of defining
relations is the following:\\

$\sigma_i \sigma_{i+1} \sigma_i = \sigma_{i+1} \sigma_i \sigma_{i+1}, ~~i = 1, 2, \ldots, n-2;$\\

$\sigma_i \sigma_{j} = \sigma_{j} \sigma_i, ~~| i - j | > 1;$\\

$ \rho \sigma_i  = \sigma_i \rho, ~~i \not= 1;$\\

$\sigma_1^{-1} \rho \sigma_{1}^{-1} \rho = \rho \sigma_{1}^{-1} \rho \sigma_{1};$\\

$\rho^2 = \sigma_1 \sigma_{2} \ldots \sigma_{n-2}  \sigma_{n-1}^2 \sigma_{n-2} \ldots \sigma_2 \sigma_{1}.$\\
\end{thm}
\begin{rem}
Geometrically element $\rho$ can be depicted similar to that of Figure~10 from \cite{Bellingeri}.
\end{rem}
\begin{proof} We start with the presentation of Van Buskirk \cite[p.~83]{VB} also studied in \cite{Gon-Gu}. It has the $2n-1$ generators
$
 \sigma_1, \sigma_2, \ldots, \sigma_{n-1}, \rho_1, \dots,  \rho_n,
$
subject to the following relations:

\smallskip

$(i) \ \sigma_i \sigma_{i+1} \sigma_i = \sigma_{i+1} \sigma_i \sigma_{i+1}, ~~i = 1, 2, \ldots, n-2;$\\

$(ii) \ \sigma_i \sigma_{j} = \sigma_{j} \sigma_i, ~~| i - j | > 1;$\\

$(iii) \ \rho_j \sigma_i  = \sigma_i \rho_j, ~~j \not= i, i+1;$\\

$(iv) \ \rho_i = \sigma_{i} \rho_{i+1} \sigma_{i};$\\

$(v) \ \rho_{i+1}^{-1} \rho_i^{-1} \rho_{i+1} \rho_i = \sigma_{i}^{2} ;$\\

$(vi) \ \rho_1^2 = \sigma_1 \sigma_{2} \ldots \sigma_{n-2}  \sigma_{n-1}^2 \sigma_{n-2} \ldots \sigma_2 \sigma_{1}.$\\
\smallskip

Let us show at first that 
the system (i) - (vi) is equivalent
 to the system (i) - (iv), (vi) and the following relations
\begin{equation}
\sigma_i^{-1} \rho_i \sigma_{i}^{-1} \rho_i = \rho_i \sigma_{i}^{-1} \rho_i \sigma_{i},
\ i=1, \dots, n-1.
\label{eq:8}
\end{equation}
We multiply the equality (\ref{eq:8}) by 
$\sigma_i\rho_i^{-1}\sigma_i\rho_i^{-1}$ 
on the left-hand side and we obtain
\begin{equation}
\sigma_i\rho_i^{-1}\sigma_i\rho_i^{-1}\sigma_i^{-1} \rho_i \sigma_{i}^{-1} \rho_i =  \sigma_{i}^2,
\ i=1, \dots, n-1.
\label{eq:8}
\end{equation}
Then we use the expression 
 $$\rho_{i+1}= \sigma_i^{-1} \rho_i \sigma_{i}^{-1}
 $$ 
 from (iv) and we obtain (v).
Now we show by induction that we can eliminate all the equalities in (\ref{eq:8})
except the first one, i.e. for $i=1$:
\begin{equation}
\sigma_1^{-1} \rho_1 \sigma_{1}^{-1} \rho_1 = \rho_1 \sigma_{1}^{-1} \rho_1 \sigma_{1}
\label{eq:7}
\end{equation}
In other words that relations (\ref{eq:8}) for $i=2, \dots,$
$ n-1$
are  consequences of relations  (ii) - (iv) and (\ref{eq:7}).
For $i=2$  we start with
(\ref{eq:7}) and multiply it by $\sigma_1^{-1}\sigma_2^{-1}$ 
on the left-hand side  and by
$\sigma_2^{-1}\sigma_1^{-1}$ on the right-hand side, we get
\begin{equation*}
\sigma_1^{-1}\sigma_2^{-1}\sigma_1^{-1} \rho_1 \sigma_{1}^{-1} \rho_1
\sigma_2^{-1}\sigma_1^{-1}=
\sigma_1^{-1}\sigma_2^{-1}\rho_1 \sigma_{1}^{-1} \rho_1 \sigma_{1}
\sigma_2^{-1}\sigma_1^{-1}.
\end{equation*}
We apply relations (ii) to this relation on the right-hand side 
 and on the left-hand side, we obtain
\begin{equation*}
\sigma_2^{-1}\sigma_1^{-1}\sigma_2^{-1} \rho_1 \sigma_{1}^{-1} \rho_1
\sigma_2^{-1}\sigma_1^{-1}=
\sigma_1^{-1}\sigma_2^{-1}\rho_1 \sigma_{1}^{-1} \rho_1 \sigma_{2}^{-1}
\sigma_1^{-1}\sigma_2.
\end{equation*}
Further we apply relation (iii) to permute $\rho_1$ and $\sigma_2^{-1}$ in all
four appearances of $\rho_1$ in the last relation, we get
\begin{equation*}
\sigma_2^{-1}\sigma_1^{-1}\rho_1 \sigma_2^{-1} \sigma_{1}^{-1}
\sigma_2^{-1}\rho_1 \sigma_1^{-1}=
\sigma_1^{-1}\rho_1 \sigma_2^{-1} \sigma_{1}^{-1} \sigma_{2}^{-1}
\rho_1 \sigma_1^{-1}\sigma_2.
\end{equation*}
Apply now  relation (ii) to the middle parts of both sides of the last relation,
and obtain
\begin{equation*}
\sigma_2^{-1}\sigma_1^{-1}\rho_1 \sigma_1^{-1} \sigma_{2}^{-1}
\sigma_1^{-1}\rho_1 \sigma_1^{-1}=
\sigma_1^{-1}\rho_1 \sigma_1^{-1} \sigma_{2}^{-1} \sigma_{1}^{-1}
\rho_1 \sigma_1^{-1}\sigma_2.
\end{equation*}
Use now relation (iv): $\rho_2 = \sigma_{1}^{-1} \rho_{1} \sigma_{1}^{-1}$
 and obtain
\begin{equation*}
\sigma_2^{-1} \rho_2 \sigma_{2}^{-1} \rho_2 = \rho_2 \sigma_{2}^{-1} \rho_2 \sigma_{2}.
\end{equation*}
This is relation (\ref{eq:8}) for $i=2$. Suppose now that for $i$ our statement is
true: the relation 
\begin{equation*}
\sigma_i^{-1} \rho_i \sigma_{i}^{-1} \rho_i = \rho_i \sigma_{i}^{-1} \rho_i \sigma_{i}.
\end{equation*}
is a consequence of relations  (ii) - (iv) and (\ref{eq:7}).
Multiplying this relation by $\sigma_i^{-1}\sigma_{i+1}^{-1}$ 
on  the left-hand side and by
$\sigma_{i+1}^{-1}\sigma_i^{-1}$ on the right-hand side and applying relations
(ii) - (iv) as before we obtain relation (\ref{eq:8}) for $i+1$. So all relations (v)
can be replaced by one relation (\ref{eq:7}).

Let us consider now relations (iii) and show that all of them are consequences of relations (i), (ii), (iv) and relations
\begin{equation}
\rho_1 \sigma_i  = \sigma_i \rho_1, ~~i \not= 1.
\label{eq:9}
\end{equation}
Really, let $j>1$, then it follows from (iv) that
\begin{equation*}
\rho_j   = \sigma_{j-1}^{-1}\sigma_{j-2}^{-1}\dots \sigma_1^{-1} \rho_1\sigma_1^{-1}
\sigma_{j-2}^{-1}\sigma_{j-1}^{-1}.
\end{equation*}
Consider $\sigma_i\rho_j$. Let $i<j-1$, then using relations (i), (ii) and
(\ref{eq:9}) we have
\begin{multline*}
\sigma_i\rho_j   = \sigma_i\sigma_{j-1}^{-1}\sigma_{j-2}^{-1}\dots \sigma_1^{-1} \rho_1\sigma_1^{-1}\dots \sigma_{j-2}^{-1}\sigma_{j-1}^{-1}= \\
\sigma_{j-1}^{-1}\sigma_{j-2}^{-1}\dots \sigma_1^{-1}\sigma_{i+1}
\rho_1\sigma_1^{-1}\dots \sigma_{j-2}^{-1}\sigma_{j-1}^{-1}= \\
\sigma_{j-1}^{-1}\sigma_{j-2}^{-1}\dots \sigma_1^{-1}\rho_1\sigma_{i+1}
\sigma_1^{-1}\dots \sigma_{j-2}^{-1}\sigma_{j-1}^{-1}=\\
\sigma_{j-1}^{-1}\sigma_{j-2}^{-1}\dots \sigma_1^{-1}\rho_1
\sigma_1^{-1}\dots \sigma_{j-2}^{-1}\sigma_{j-1}^{-1}\sigma_{i}= \rho_j\sigma_i
\end{multline*}
If $i>j$, then using relations (i), and
(\ref{eq:9}) we have
\begin{multline*}
\sigma_i\rho_j   = \sigma_i\sigma_{j-1}^{-1}\sigma_{j-2}^{-1}\dots \sigma_1^{-1} \rho_1\sigma_1^{-1}\dots \sigma_{j-2}^{-1}\sigma_{j-1}^{-1}= \\
\sigma_{j-1}^{-1}\sigma_{j-2}^{-1}\dots \sigma_1^{-1}\sigma_{i}
\rho_1\sigma_1^{-1}\dots \sigma_{j-2}^{-1}\sigma_{j-1}^{-1}= \\
\sigma_{j-1}^{-1}\sigma_{j-2}^{-1}\dots \sigma_1^{-1}\rho_1\sigma_{i}
\sigma_1^{-1}\dots \sigma_{j-2}^{-1}\sigma_{j-1}^{-1}=\\
\sigma_{j-1}^{-1}\sigma_{j-2}^{-1}\dots \sigma_1^{-1}\rho_1
\sigma_1^{-1}\dots \sigma_{j-2}^{-1}\sigma_{j-1}^{-1}\sigma_{i}= \rho_j\sigma_i
\end{multline*}
Hence all relations (iii) are consequences of relations (i), (ii), (iv) and
(\ref{eq:9}). So, we can delete generators
$
  \rho_2, \dots,  \rho_n,
$
and relations (iv) from the presentation and replace relations (iii) and (v) by
relations (\ref{eq:9}) and  (\ref{eq:7}) respectively.
\end{proof}
There is a canonical homomorphism $\tau : B_n(\RP^2) \longrightarrow \Sigma_n,$ $\tau(\sigma_i) = (i,i+1),$ $\tau(\rho) = e.$ The kernel $\mathrm{Ker} (\tau)$ is the pure
braid group $P_n(\RP^2).$ This group was studied in \cite{Gon-Gu}.
We will find a presentation for $P_3(\RP^2)$ which we shall use later.
Consider at first  the group $B_2(\RP^2)$. We have
$$
B_2(\RP^2) = \langle \rho, \sigma_1 ~|~\sigma_1^{-1} \rho \sigma_{1}^{-1} \rho = \rho \sigma_{1}^{-1} \rho \sigma_{1}, ~~~\rho^2 = \sigma_1^2 \rangle.
$$
This group has order 16  and $ P_2(\RP^2)$
is isomorphic to the quaternion group ${\bold Q}_8$ of order 8 \cite{VB}.
Relation $\rho^2 = \sigma_1^2$ gives that the corresponding
pure braid group $ P_2(\RP^2)$ is normally generated by $\rho$
and makes it possible not to use the canonical generator of the
pure braid group $A_{12} = \sigma_1^2$.
Let us define the following element of $ P_2(\RP^2)$:
$$
u = \sigma_1 \rho \sigma_1^{-1}.
$$
Reidemeister method gives
 the following presentation:
\begin{equation}
P_2(\RP^2) = \langle \rho, u ~|~\rho u \rho = u, ~~~\rho^{2} = u^{2} \rangle.
\label{eq:presRP2}
\end{equation}

Consider now the case $n = 3.$ We have:
$$
B_3(\RP^2) = \langle \rho, \sigma_1, \sigma_2~|~\sigma_1 \sigma_{2} \sigma_1 = \sigma_{2} \sigma_1 \sigma_{2}, ~~~\rho \sigma_{2} = \sigma_{2} \rho,
$$
$$
\sigma_1^{-1} \rho \sigma_{1}^{-1} \rho = \rho \sigma_{1}^{-1} \rho \sigma_{1}, ~~~\rho^2 = \sigma_1 \sigma_{2}^2 \sigma_{1} \rangle.
$$
To construct a presentation for $P_3(\RP^2)$ we  use the Reidemeister-Schreier method.
As representatives of cosets of the normal subgroup  $P_3(\RP^2)$ in the group
$B_3(\RP^2)$ we take the elements: $e$, $\sigma_1$, $\sigma_2$, $\sigma_2\sigma_1$
$\sigma_1\sigma_2$, $\sigma_1\sigma_2\sigma_1$. Then by \cite[Theorem~2.7]{MKS}
the group $P_3(\RP^2)$ is generated by
elements
\begin{equation*}
ka\overline{(ka)}^{-1},
\end{equation*}
where $a\in \{\rho, \sigma_1, \sigma_2\}$, $k\in \{$
$e$, $\sigma_1$, $\sigma_2$, $\sigma_2\sigma_1$,
$\sigma_1\sigma_2$, $\sigma_1\sigma_2\sigma_1$ $\}$ and the bar denotes the mapping from words to their coset representatives \cite[p.~88]{MKS}.
Having in mind that $\sigma_2\rho\sigma_2^{-1}= \rho$,
we obtain that the group $P_3(\RP^2)$ is generated by
$$
\rho,~~u = \sigma_1 \rho \sigma_1^{-1},~~w = \sigma_2 \sigma_1 \rho \sigma_1^{-1} \sigma_2^{-1},~~ A_{12}, ~~A_{23} = \sigma_2^2, ~~
A_{13} = \sigma_2 \sigma_1^2 \sigma_2^{-1}.
$$
The following set of defining relations is obtained by application of
Reidemeister-Schreier method \cite[Theorem~2.9]{MKS}:

$$
A_{12} A_{13} A_{12}^{-1} = A_{23}^{-1} A_{13} A_{23},~~~A_{12} \left(A_{13} A_{23} \right) A_{12}^{-1} = A_{13} A_{23},
$$

$$
\rho A_{23} \rho^{-1} = A_{23},~~~ u \left( A_{23}^{-1} A_{13} A_{23} \right) u^{-1} = A_{23}^{-1} A_{13} A_{23},
$$

\begin{equation}
\rho (A_{13}^{-1} w^{-1} A_{13} ) \rho^{-1} = w^{-1} A_{13},~~~ \rho \left( A_{13}^{-1} w \right) \rho^{-1} = w,
\label{eq:conroA13w}
\end{equation}

$$
\rho (A_{12}^{-1} u) \rho^{-1} = u,
$$
\begin{equation}
 ~~~u ( A_{23}^{-1} w^{-1} A_{23} ) u^{-1} = w^{-1} A_{23},~~~ u \left( A_{23}^{-1} w \right) u^{-1} = w,
 \label{eq:conuAw}
 \end{equation}

\begin{equation}
A_{23}^{-1} A_{13} A_{23} A_{12} = \rho^2,~~~ A_{12} A_{13} = \rho^2,
~~~A_{12} A_{23} = u^2,
\label{eq:square}
\end{equation}
\begin{equation}
~~~ A_{13} A_{23} = w^2.
\label{eq:A13A23}
\end{equation}
From these relations we have the following formulas for conjugation by $A_{12}$, $\rho$, $u$:
$$
A_{12} A_{13} A_{12}^{-1} = A_{23}^{-1} A_{13} A_{23},~~~A_{12}  A_{23}  A_{12}^{-1} = A_{23}^{-1} A_{13}^{-1} A_{23} A_{13} A_{23},$$

\begin{equation}
A_{12} w A_{12}^{-1} = w,
\label{eq:conA12w}
\end{equation}

\begin{equation}
\rho w \rho^{-1} = w^{-1} A_{13}^{-1} w^2,~~~\rho A_{13}
\rho^{-1} = w^{-1} A_{13}^{-1} w,
\label{eq:conrowA}
\end{equation}

$$
~~~\rho A_{23} \rho^{-1} = A_{23},
$$

\begin{equation}
u w u^{-1} = w^{-1} A_{23}^{-1} w^2,~~~
u A_{23} u^{-1} = w^{-1} A_{23}^{-1} w,
\label{eq:conuw}
\end{equation}
\begin{equation}
u A_{13} u^{-1} = w^{-1}A_{23}^{-1}wA_{23}^{-1}wA_{23}w.
\label{eq:conuA13}
\end{equation}

\begin{rem} Relation (\ref{eq:conA12w}) can be more easily seen directly from
relations of $B_3(\RP^2)$.  Relations (\ref{eq:conrowA}) are obtained from
relations (\ref{eq:conroA13w}). Relations (\ref{eq:conuw}) are obtained from
relations (\ref{eq:conuAw}).  Relation (\ref{eq:conuA13}) is obtained from
relations (\ref{eq:A13A23}) and (\ref{eq:conuw}).
\end{rem}

We see from these formulas that the subgroup
$$
U_3(\RP^2) = \langle w, A_{13}, A_{23} ~||~ A_{13} A_{23} = w^2 \rangle
$$ is normal in $P_3(\RP^2)$.
Geometrically it can be identified  with $\pi_1(\RP^2\smallsetminus\{p_1,p_2\})$
which is included in short exact sequence, see diagram (\ref{equation3.2}):
\begin{equation*}
\pi_1(\RP^2\smallsetminus\{p_1,p_2\})\rInto^{{i_3}_*}
P_3(\RP^2)\rOnto^{d_3} P_2(\RP^2)
\end{equation*}
 and so $U_3(\RP^2)$ is the free group of rank 2 and $P_3(\RP^2) / U_3(\RP^2) \simeq P_2(\RP^2)$.

We can exclude the generators $A_{12}$, $A_{13}$ from the list of generators $P_3(\RP^2)$, using the formulas:
\begin{equation}
A_{12} = u \rho^{-1} u^{-1} \rho,
\label{eq:A12}
\end{equation}
$$
~~~A_{13} = w^2 A_{23}^{-1}.
$$
Inserting this formulas in the rest of relations of the group $P_3(\RP^2)$  we
obtain the following statement.
\begin{lem} \label{p3rp2}
The group $P_3(\RP^2)$ can be generated by elements
$$
\rho, u, w, A_{23}
$$
and has the following relations:

$$
\rho w \rho^{-1} = w^{-1} A_{23},~~~\rho A_{23} \rho^{-1} = A_{23}, \eqno{(1)}
$$

$$
\rho^{-1} w \rho = A_{23} w^{-1},~~~\rho^{-1} A_{23} \rho = A_{23}, \eqno{(1')}
$$

$$
u w u^{-1} = w^{-1} A_{23}^{-1} w^2,~~~u A_{23} u^{-1} = w^{-1} A_{23}^{-1} w, \eqno{(2)}
$$

$$
u^{-1} w u = A_{23}^{-1} w,~~~u^{-1} A_{23} u = A_{23}^{-1} w A_{23}^{-1} w^{-1} A_{23}, \eqno{(2')}
$$

$$
\rho^{-1} u \rho^{-1}  u^{-1} = w A_{23}^{-1} w,~~~u^{-1} \rho^{-1} u^{-1} \rho = A_{23}^{-1}. \eqno{(3)}
$$
\end{lem}
\begin{proof}
The first relation in (3) follows from (\ref{eq:A12}), the second relation in 
(\ref{eq:square}) and the second relation in  (\ref{eq:conrowA}). 
The second relation in (3) follows from (\ref{eq:A12})and  the third relation in 
(\ref{eq:square}). 
To prove that the statement of Lemma gives a presentation 
of $P_3(\RP^2)$ denote by
$P$ the group which presentation is given by these generators and relations.
There exists an evident homomorphism  
$$\phi: P\to P_3(\RP^2).$$
  The subgroup $U_3(\RP^2)$ generated by $w, A_{2,3}$ is a free subgroup in
  $P$ as it is free after the mapping by $\phi$. The quotient $P/U_3(\RP^2)$
  is isomorphic to $P_2(\RP^2)$ (relations (3)); so $\phi$ becomes an
isomorphism after comparison of exact sequences:
\begin{equation*}
\begin{diagram}
U_3(\RP^2)&\rInto^{} &P&\rOnto^{}& P_2(\RP^2)\\
\dTo>{}&&\dTo>{\phi}&&\dTo>{}\\
U_3(\RP^2&\rInto^{{i_3}_*}&P_3(\RP^2))&\rOnto^{d_2}& P_2(\RP^2).\\
\end{diagram}
\end{equation*}  
\end{proof}
\begin{rem}
1) Certainly, relations $(1')$ and $(2')$ follow from relations $(1)$  and $(2)$ respectively.

2) Similar presentation was constructed in \cite[p.~765]{Gon-Gu} but in the list of relations there, in  the forth relation of formula  $(3)$ instead of
$$
\rho_2^{-1} B_{2,3} \rho_2 = B_{2,3}^{-1} \rho_3 B_{2,3} \rho_3^{-1} B_{2,3}
$$
it should  be
$$
\rho_2^{-1} B_{2,3} \rho_2 = B_{2,3}^{-1} \rho_3 B_{2,3}^{-1} \rho_3^{-1} B_{2,3}.
$$
\end{rem}

Let us introduce new generators $a = \rho w,$ $b = w u$ then from Lemma~\ref{p3rp2}
 we get
the following statement.
\begin{lem}
The group $P_3(\RP^2)$ can be generated by elements
$$
a, b, w, A_{23}
$$
and has the following relations:

$$
a w a^{-1} = w^{-1} A_{23},~~~a A_{23} a^{-1} = w^{-1} A_{23} w, \eqno{(4)}
$$

$$
a^{-1} w a =  w^{-1} A_{23},~~~a^{-1} A_{23} a = w^{-1} A_{23} w, \eqno{(4')}
$$

$$
b w b^{-1} = A_{23}^{-1} w,~~~b A_{23} b^{-1} = A_{23}^{-1}, \eqno{(5)}
$$

$$
b^{-1} w b = A_{23}^{-1} w,~~~b^{-1} A_{23} b = A_{23}^{-1}, \eqno{(5')}
$$

$$
b a b^{-1}  = a^{-1},~~~a^{2} = b^{2}. \eqno{(6)}
$$
In particular, $\langle a, b \rangle \simeq P_2(\RP^2) \leq P_3(\RP^2).$ $\square$
\end{lem}

From this Lemma we have the following statement.
\begin{prop}
There exists the following splitting short exact sequence
$$
1 \longrightarrow U_3(\RP^2) \longrightarrow P_3(\RP^2)
\stackrel{d_3} {\longrightarrow} P_2(\RP^2) \longrightarrow 1
$$
and $P_3(\RP^2) = U_3(\RP^2) \leftthreetimes P_2(\RP^2)$.
\label{prop:split}
\end{prop}
\hfill $\Box$

\noindent
This proposition was proved
in \cite{Gon-Gu}. It was proved also there that for $n = 2, 3$ and for all $m \geq 4$ the short exact sequence:
$$
1 \longrightarrow P_{m-n}(\RP^2 \setminus \{ x_1, x_2, \ldots, x_n \}) \longrightarrow P_m(\RP^2) \longrightarrow P_n(\RP^2) \longrightarrow 1
$$
does not split.

\subsection{$3$-strand Brunnian braids on the Projective Plane}

In order to pass to the brunnian braids recall the geometric interpretations for the generators $\rho$, $u$, $w$. We represent $\RP^2$ as 2-gon $L$ where opposite points on two edges are
identified by standard manner. In the braid $\rho$ the second and the third
strings are just two parallel lines. It's  first strand which goes through the edge of L. The braids $u$ and
$w$ are defined by similar manner. In  $u$ the second strand passes through the edge
and in $w$ the third one.
The braid $A_{23}$ is defined as in the braid group of disk.
Remind that the presentation of $P_2(\RP^2)$ is given by the formula
(\ref{eq:presRP2}).
Hence the maps
$$
d_1, d_2, d_3 : P_3(\RP^2) \longrightarrow P_2(\RP^2)
$$
act on the generators by the rules:
$$
d_1 :
\left\{%
\begin{array}{lcl}
a & \longrightarrow & u, \\
b & \longrightarrow & u\rho, \\
A_{23} & \longrightarrow & A_{12}, \\
w & \longrightarrow & u, \\
\end{array}%
\right.~~~~
d_2 :
\left\{%
\begin{array}{lcl}
a & \longrightarrow & \rho u, \\
b & \longrightarrow &  u, \\
A_{23} & \longrightarrow & 1, \\
w & \longrightarrow & u, \\
\end{array}%
\right.~~~~
d_3 :
\left\{%
\begin{array}{lcl}
a & \longrightarrow & \rho, \\
b & \longrightarrow & u, \\
A_{23} & \longrightarrow & 1, \\
w & \longrightarrow & 1. \\
\end{array}%
\right.~~~~
$$

From the exact sequence of Proposition~\ref{prop:split} we see that
$\Brun_3(\RP^2)$ is a subgroup of $U_3(\RP^2)$:
$$
\Brun_3(\RP^2) \leq U_3(\RP^2),
$$
so in our study of Brunnian braids on $\RP^2$ we can restrict ourselves to
look at $U_3(\RP^2)$ and the action of $d_1$ and $d_2$ on it. We write
the action of $d_3$ as a supplementary information.

We have
$$
d_1(w^4) = d_2(w^4) = u^4,~~~d_3(w^4) = 1
$$
and since $u^4 = 1$ in $P_2(\RP^2)$ then $w^4 \in \Brun_3(\RP^2).$ Similarly
$$
d_1(A_{23}^2) = A_{12}^2,~~~d_2(A_{23}^2) = d_3(A_{23}^2) = 1
$$
and since $A_{12}^2 = \sigma_1^4 = 1$ in $P_2(\RP^2)$
(se formula (\ref{eq:presRP2})) then $A_{23}^2 \in \Brun_3(\RP^2).$
For the commutator $[w, A_{23}]$ we have
$$
d_1([w, A_{23}]) = [u, A_{12}],~~~d_2([w, A_{23}]) = d_3([w, A_{23}]) = 1
$$
and since $A_{12}$ lies in the center of $P_2(\RP^2)$ then $d_1([w, A_{23}]) = 1$ and $[w, A_{23}] \in \Brun_3(\RP^2).$

Now we are going to determine a free basis for $\Brun_3(\RP^2)$.
\begin{lem}
Let $F(S)$ be the free group (freely) generated by the set $S$. Given $x\in S$, let $C_q(x)\cong \Z/q$ be the cyclic group of order $q$ generated by $x$ and let $p_x\colon F(S)\to C_q(x)$ be the group homomorphism with $p(y)=1$ for $y\not=x\in S$ and $p_x(x)=x$. Then $\Ker(p_x)$ has a free basis
$$\{x^q, y, [y,x^j] \ | \, y\in S,\ y\not=x \ 1\leq j\leq q-1\}.$$
\end{lem}
\begin{proof}
By using Schreier system, $\Ker(p_x)$ has a free basis $$\{x^q, x^{-j}yx^j \ | \ y\in S, \ y\not=x, 0\leq j\leq q-1\}$$ which is equivalent to the generating set in the statement as $$[y,x^j]=y^{-1}(x^{-j}yx^j)$$ and hence the assertion.
\end{proof}

\begin{prop}\label{brun3-rp2}
As a subgroup of $B_3(\RP^2)$, $\Brun_3(\RP^2)$ has a free basis given by
$$
\begin{array}{ccc}
 x_2^2,&x_1^4, & [x_1^4,x_2], \\
 &[x_2,x_1],& [[x_2,x_1],x_2],\\
&[x_2,x_1^2], &[[x_2,x_1^2],x_2],\\
&[x_2,x_1^3], &[[x_2,x_1^3],x_2],\\
\end{array}
$$
where $x_1=w$ and $x_2=A_{2,3}$.
\end{prop}
\begin{proof}
Consider the projection $p_{x_1}\colon F(x_1,x_2)\to C_4(x_1)$. 
(It is $d_2$ in our case.)
By the above lemma, $\Ker(p_{x_1})$ has a free basis given by
$$
S=\{x_1^4, x_2, [x_2,x_1], [x_2,x_1^2], [x_2,x_1^3]\}
$$
The assertion follows by applying the above lemma to the projection $p_{x_2}\colon F(x_1,x_2)\to C_2(x_2)$ 
($d_1$ in our case) restricted to the subgroup $F(S)=\Ker(p_{x_1})$.
\end{proof}

Let us describe the quotient groups $P_3(\RP^2) / \Brun_3(\RP^2)$ and \newline
$B_3(\RP^2) / \Brun_3(\RP^2)$.

\begin{prop}
1)
Let $\overline{w}$ and $\overline{A}$ be the images of $w$ and $A_{23}$ respectively after natural projection
$$
U_3(\RP^2) \longrightarrow U_3(\RP^2) / \Brun_3(\RP^2)
$$
then
$$
U_3(\RP^2) / \Brun_3(\RP^2) = \langle \overline{w}, \overline{A}~||~ \overline{w}^{4} = \overline{A}^2 = 1,~~~\overline{A} \overline{w} =
\overline{w} \overline{A} \rangle \simeq \mathbb{Z}_4 \oplus \mathbb{Z}_2.
$$

2)
The quotient $P_3(\RP^2) / \Brun_3(\RP^2)$ has order 64 and is the semi-direct product
$$
P_3(\RP^2) / \Brun_3(\RP^2) = (U_3(\RP^2) / \Brun_3(\RP^2)) \leftthreetimes P_2(\RP^2).
$$
More precisely  $P_3(\RP^2) / \Brun_3(\RP^2)$ is generated by
$$
 \overline{w}, \overline{A}, a, b
$$
and has defining relations:
\begin{equation*}
 \overline{w}^{4} = \overline{A}^2 = 1,~~~\overline{A} \overline{w} =
\overline{w} \overline{A},~~~b a b^{-1} = a^{-1},~~~a^2 = b^2,
\end{equation*}
\begin{equation*}
a^{-1} \overline{w} a = \overline{w}^{-1} \overline{A},~~~a^{-1} \overline{A} a = \overline{A}, \ \eqno{(1)}
\end{equation*}
\begin{equation*}
a \overline{w} a^{-1} = \overline{w}^{-1} \overline{A},~~~a \overline{A} a^{-1} = \overline{A}, \ \eqno{(1^\prime)}
\end{equation*}
\begin{equation*}
b^{-1} \overline{w} b = \overline{w} \overline{A},~~~b^{-1} \overline{A} b = \overline{A}, \eqno{(2)}
\end{equation*} 
\begin{equation*}
b \overline{w} b^{-1} = \overline{w} \overline{A},~~~b \overline{A} b^{-1} = \overline{A}. \eqno{(2')}
\end{equation*}
\hfill $\Box$
\end{prop}
\begin{rem}
Relations with primes are equivalent to those without prime.
\end{rem} 

Using the short exact sequence
\begin{equation}
1 \longrightarrow P_3(\RP^2) \longrightarrow B_3(\RP^2) \longrightarrow \Sigma_3 \longrightarrow 1.
\label{eq:pbsig}
\end{equation}
we want to describe $B_3(\RP^2)$ as extension of $P_3(\RP^2)$ by $\Sigma_3$.
\begin{prop}
 The group $B_3(\RP^2)$ can be presented as having the generators
$$
a,~~b,~~w,~~A_{23},~~\sigma_1,~~\sigma_2,
$$
and the following relations:
\begin{equation}
\sigma_1^2 = a^2 w^{-2},~~~\sigma_2^2 = A_{23}.
\label{eq:00}
\end{equation}
\begin{equation}
\sigma_1^{-1} a \sigma_1 = b A_{23}^{-1}, \label{eq:1.1}
\end{equation}
\begin{equation}
\sigma_1^{-1} b
\sigma_1 = a w^{-1} A_{23} w^{-1}, \label{eq:1.2}
\end{equation}
\begin{equation}
\sigma_1^{-1} w \sigma_1 = w, \label{eq:1.3}
\end{equation}
\begin{equation}
\sigma_1^{-1} A_{23} \sigma_1 = w^2 A_{23}^{-1}, \label{eq:1.4}
\end{equation}
\begin{equation}
 \sigma_2^{-1} a
\sigma_2 = a b (w^{-1} A_{23})^{2}, \label{eq:2.1}
\end{equation}
\begin{equation}
\sigma_2^{-1} b \sigma_2 = b A_{23}, \label{eq:2.2}
\end{equation}
\begin{equation}
\sigma_2^{-1} w \sigma_2 = b w^{-1} A_{23}, \label{eq:2.3}
\end{equation}
\begin{equation}
 \sigma_2^{-1} A_{23} \sigma_2 = A_{23}. \label{eq:2.4}
\end{equation}
\end{prop}
\begin{proof}
The first relation in (\ref{eq:00}) follows from the definition of 
elements $a$ and $w$ and relations of the presentation of $B_3(\RP^2)$ 
with generators $\rho, \sigma_1$ and $\sigma_2$. The second relation
in (\ref{eq:00}) is the definition of $A_{23}$.

To construct the formulas of conjugation we can take the corresponding 
relations from the paper of Van Buskirk \cite{VB} and rewrite them in our generators of $P_3(\RP^2)$.
On the other side, when the formulas are written  one can prove them. Let us do it.
At first let us prove (\ref{eq:1.4}). We start by the two equal expressions of $A_{13}$:
\begin{equation}
\sigma_1^{-1} \sigma_2^2 \sigma_1 = \sigma_2 \sigma_1^2  \sigma_{2}^{-1},
\label{eq:2A23}
\end{equation}
this is true in  $B_3(\RP^2)$. We insert $\sigma_2\sigma_2^{-1}$ in the right hand part of (\ref{eq:2A23}): 
\begin{equation*}
\sigma_1^{-1} \sigma_2^2 \sigma_1 = 
\sigma_2 \sigma_1^2 \sigma_2 \sigma_{2}^{-2},
\end{equation*}
Then we use the relation $\rho^2 = \sigma_1 \sigma_2^2 \sigma_1$ 
from the presentation of $B_3(\RP^2)$: 
$$
\sigma_1^{-1} \sigma_2^2 \sigma_1 = \sigma_2 \sigma_1 \rho^2 \sigma_1^{-1} \sigma_2^{-1} A_{23}^{-1}.
$$
Since $A_{23} = \sigma_2^2$ and $w = \sigma_2 \sigma_1 \rho \sigma_1^{-1} \sigma_2^{-1}$  we have
\begin{equation*}
\sigma_1^{-1} A_{23} \sigma_1 = w^2 A_{23}^{-1}, 
\end{equation*}

Relation (\ref{eq:1.3}). We start with the definition of $w$:
$$
\sigma_2 \sigma_1  \rho \sigma_1^{-1} \sigma_2^{-1} = w.
$$
Since $\rho \sigma_2 = \sigma_2 \rho$ we have 
\begin{multline*}
\sigma_2 \sigma_1  \rho \sigma_1^{-1} \sigma_2^{-1} =
(\sigma_2 \sigma_1 \sigma_2^{-1}) \rho (\sigma_2 \sigma_1^{-1} \sigma_2^{-1}) =\\(\sigma_1^{-1} \sigma_2 \sigma_1) \rho (\sigma_1^{-1} \sigma_2^{-1} \sigma_1) =\sigma_1^{-1} w \sigma_1, 
\end{multline*}
and so:
$$
\sigma_1^{-1} w \sigma_1 = w.
$$

Relation (\ref{eq:1.2}). We start with relation ($1'$) from Lemma\ref{p3rp2}
$$
\rho^{-1} w \rho =  A_{23} w^{-1},
$$
which is equivalent to
$$
w \rho = \rho  A_{23} w^{-1}.
$$
Since $\sigma_1^{-1} w \sigma_1 = w$ and $\sigma_1^{-1} u \sigma_1 = \rho$ 
we have
$$
\sigma_1^{-1} w u \sigma_1 = (\rho w) w^{-1} A_{23} w^{-1}.
$$
Using the definition of $a$ and $b$ we obtain (\ref{eq:1.2}).

Relation (\ref{eq:1.1}). We start with equality
$$
 w = (A_{23} w^{-1}) (w u A_{23}^{-1} w)
$$
and apply  the conjugation formulas ($1^\prime$) and (3) 
($1'$) from Lemma\ref{p3rp2}, we have
$$
 w = (\rho^{-1} w \rho) (\rho^{-1} u \rho^{-1} u^{-1}).
$$
what is equivalent to
\begin{equation}
 w = \rho^{-1} w u \rho^{-1} u^{-1}.
 \label{wrhou}
\end{equation}
We rewrite the first equation in  (3)  from Lemma\ref{p3rp2}
in the form
$$
\rho (w A_{23}^{-1} w ) u \rho u^{-1}=1
$$
and multiply by $\rho (w A_{23}^{-1} w ) u \rho u^{-1}$
the right hand side of (\ref{wrhou}), we obtain
$$
 w = \rho^{-1} w u \rho^{-1} u^{-1} \rho (w A_{23}^{-1} w ) u \rho u^{-1}.
$$
We apply (1) of Lemma\ref{p3rp2} and we get
$$
 w = \rho^{-1} w u \rho^{-1} u^{-1} \rho (\rho A_{23} w^{-2} \rho^{-1}) u \rho u^{-1}.
$$
Using the formulas
$$
A_{12} = u \rho^{-1} u^{-1} \rho,~~~A_{13}^{-1} = A_{23} w^{-2}.
$$
we obtain
$$
 w = \rho^{-1} w A_{12} \rho A_{13}^{-1} A_{12}^{-1}
$$
or
$$
\rho w = w A_{12} \rho A_{13}^{-1} A_{12}^{-1}.
$$
Conjugating it by $\sigma_1^{-1}$ we have
$$
\sigma_1^{-1} (\rho w) \sigma_1 = w u A_{23}^{-1}.
$$
and this is (\ref{eq:1.1}).

Formula (\ref{eq:2.4}) follows from $A_{23} = \sigma_2^2$.

Relation (\ref{eq:2.3})
We start with the first relation in ($2'$) of Lemma\ref{p3rp2}
and we rewrite it in equivalent forms
$$
u^{-1} w u = A_{23}^{-1} w \Leftrightarrow 1 = u^{-1} w u w^{-1} A_{23} \Leftrightarrow  u = w u w^{-1} A_{23} 
$$
or
$$
\sigma_2^{-1} (\sigma_2 \sigma_1 \rho \sigma_1^{-1} 
\sigma_2^{-1}) \sigma_2 = w u w^{-1} A_{23}
$$
what is equivalent to  (\ref{eq:2.3}) :
$$
\sigma_2^{-1} w \sigma_2 = b w^{-1} A_{23}.
$$

Relation  (\ref{eq:2.2}). We start with the identity
$$
( b w^{-1} A_{23}) ( A_{23}^{-1} w A_{23}) = b A_{23}.
$$
Using the formula
$$
\sigma_2^{-1} u \sigma_2 = A_{23}^{-1} w A_{23},
$$
and  (\ref{eq:2.3}) we get
$$
(\sigma_2^{-1} w \sigma_2) ( \sigma_2^{-1} u \sigma_2) = b A_{23}.
$$
what is equivalent to (\ref{eq:2.2}) :
$$
\sigma_2^{-1} (w u) \sigma_2 = b A_{23}.
$$

Finally let us prove relation ((\ref{eq:2.1}).
 We start with the identity
$$
1 = (A_{23}^{-1} w )  w^{-1} A_{23}.
$$ 
Using   ($2'$) of Lemma\ref{p3rp2} we get
$$
1 = (u^{-1}  w u)  w^{-1} A_{23}
$$ 
 and then
$$
u = (  w u)  w^{-1} A_{23}.
$$ 
 We multiply this equality by $\rho w$ from the left hand side
$$
 \rho  w u  = (\rho w) (w u) w^{-1} A_{23}
$$ 
what is equivalent to
$$
 \rho  b  = a b w^{-1} A_{23}.
$$
Since $\rho$ and $\sigma_2$ commute this is the same as
$$
 (\sigma_2^{-1} \rho  \sigma_2) b  = a b w^{-1} A_{23}.
$$ 
 Multiply this equality by $w^{-1} A_{23}$ from the right hand side
$$
 (\sigma_2^{-1} \rho  \sigma_2) (b w^{-1} A_{23}) = a b (w^{-1} A_{23})^{2}.
$$ 
and use ((\ref{eq:2.3})
$$
 \sigma_2^{-1} \rho w \sigma_2 = a b (w^{-1} A_{23})^{2}.
$$
what is equivalent to ((\ref{eq:2.3}):
$$
 \sigma_2^{-1} a \sigma_2 = a b (w^{-1} A_{23})^{2}.
$$
The proof that $B_3(\RP^2)$ has a presentation as in the statement of 
the Proposition is the same as the proof of the presentation 
of Lemma\ref{p3rp2} with the help of the exact sequence (\ref{eq:pbsig}).
\end{proof}

\begin{prop}
The quotient $B_3(\RP^2) / \Brun_3(\RP^2)$ has order 384 and is an extension of $P_3(\RP^2) / \Brun_3(\RP^2)$ by $\Sigma_3$:
$$
1 \longrightarrow P_3(\RP^2) / \Brun_3(\RP^2) \longrightarrow B_3(\RP^2) / \Brun_3(\RP^2) \longrightarrow \Sigma_3 \longrightarrow 1.
$$
 $B_3(\RP^2) / \Brun_3(\RP^2)$ is generated by
$$
 \overline{w}, \overline{A}, a, b, \sigma_1, \sigma_2,
$$
and has defining relations:
$$
\sigma_1^2 = a^2 \overline{w}^{2},~~~\sigma_2^2 = \overline{A},
$$
$$
 \overline{w}^{4} = \overline{A}^2 = 1,~~~\overline{A} \overline{w} =
\overline{w} \overline{A},~~~b a b^{-1} = a^{-1},~~~a^2 = b^2,
$$
$$
a^{-1} \overline{w} a = \overline{w}^{-1} \overline{A},~~~a^{-1} \overline{A} a = \overline{A},
$$
$$
b^{-1} \overline{w} b = \overline{w} \overline{A},~~~b^{-1} \overline{A} b = \overline{A},
$$
$$
\sigma_1^{-1} a \sigma_1 = b \overline{A},~~~\sigma_1^{-1} b \sigma_1 = a \overline{A} \overline{w}^{2},~~~\sigma_1^{-1} \overline{w} \sigma_1 = \overline{w},~~~
\sigma_1^{-1} \overline{A} \sigma_1 = \overline{A} \overline{w}^{2},
$$
$$
\sigma_2^{-1} a \sigma_2 = a b \overline{w}^2,~~~\sigma_2^{-1} b \sigma_2 = b \overline{A},~~~\sigma_2^{-1} \overline{w} \sigma_2 = \overline{A} \overline{w}^{-1},~~~
\sigma_2^{-1} \overline{A} \sigma_2 = \overline{A}.
$$
\end{prop}

\section{Proof of Theorem~\ref{theorem1.2}}\label{section5}
\subsection{Some Lemmas on Free Groups}\label{section-free-group}
Let $S$ be a set and let $F(S)$ be the free group freely generated
by $S$. Let $S_0$ be a set and let $x_1,x_2,\ldots$ be additional
letters. Let $S_n=S_0\cup\{x_1,\ldots,x_n\}$ be the disjoint
union. Consider the group homomorphism
$$
d_i\colon F(S_n)\to F(S_{n-1}), \ \  1\leq i\leq n,
$$
such that
\begin{equation}\label{equation4.1}
d_i(x)=\left\{
\begin{array}{rcl}
x&\textrm{ if } & x\in S_0 \textrm{ or } x=x_j\textrm{ with } j<i,\\
1&\textrm{ if }& x=x_i,\\
x_{j-1}&\textrm{ if }& x=x_j\textrm{ with } j>i.\\
\end{array}
\right.
\end{equation}
Roughly speaking, $d_i$ is obtained by sending $x_i$ to $1$ and
keeping other generators.  The following lemma is a special case of~\cite[Theorem 4.3]{LiWu}.

\begin{lem}\label{lemma4.2}
Let $d_i\colon F(S_n)\to F(S_{n-1})$ be defined by the formula~(\ref{equation4.1}).
Then
$$
\bigcap_{j=1}^k\Ker(d_i)=[\Ker(d_1),\Ker(d_2),\ldots,\Ker(d_k)]_S
$$
for $2\leq k\leq n$.\hfill $\Box$
\end{lem}

Let $H$ be a normal subgroup of $G$. A
set $X$ of elements of $H$ is called a set of \textit{normal generators}
for $H$ in $G$ if $H$ is the normal closure of $X$ in $G$.
  We say that $H$ has \textit{finitely many normal generators} in $G$ if there is a finite set $X$ such that $H$ is the normal closure of $X$ in $G$.

\begin{lem}\label{lemma5.1}
Let $R_1$ and $R_2$ be normal subgroups of $G$. Suppose that
\begin{enumerate}
\item $R_1$ has finitely many normal generators and
\item $R_2$ has finitely many generators (in the usual sense).
\end{enumerate}
Then the commutator subgroup $[R_1,R_2]$ has finitely many normal
generators.
\end{lem}
\begin{proof}
Let $\{a_1,\ldots,a_m\}$ be a set of normal generators for $R_1$.
The set of generators for $R_1$ can be given as $\{g^{-1}a_ig \ | \
1\leq i\leq m, \ g\in G\}$. Let $\{b_1,\ldots, b_n\}$ be a set of
generators   for $R_2$. Let $H$ be the normal closure of
$$
\{ [a_i,b_j]\ | \ 1\leq i\leq m,\ 1\leq j\leq n\}.
$$
Now take any $r\in R_2$, $r=b_{i_1}\ldots b_{i_k}$.
Then
$$
[a_i,r] = [a_i, b_{i_1}] \,g_1[a_i, b_{i_1}] g_1^{-1} \ldots
g_j[a_i, b_{i_{j-1}}] g_j^{-1}
\ldots g_{k-1}[a_i, b_{i_k}] g_{k-1}^{-1},
$$
where $g_j=b_{i_1}\dots b_{i_{j}}$. So
$[a_i,r]\in H$ for any $r\in R_2$.
Now
$$
[g^{-1}a_ig, b_j]=g^{-1}[a_i,gb_jg^{-1}]g\in H
$$
because $gb_jg^{-1}\in R_2$. We get that  $[R_1,R_2]=H$.
\end{proof}

\begin{lem}\label{lemma5.2}
Let $M$ be a connected compact $2$-manifold with nonempty boundary. Let $n\geq 2$. Then the subgroup
$$
\bigcap_{i=1}^k \Ker(d_i\colon P_n(M)\to P_{n-1}(M))\cap\Ker(d_n\colon P_n(M)\to P_{n-1}(M))
$$
has finitely many normal generators in $P_n(M)$ for each $1\leq k\leq n-1$.
\end{lem}
\begin{proof}
The proof is given by induction on $k$. The assertion holds for $k=1$ by Lemma~\ref{lemma3.10}. Suppose that the assertion holds for $k-1$. Consider the short exact sequence of groups
$$
\pi_1(M\smallsetminus\{p_1,\ldots,p_{n-1}\})\rInto^{i_*} P_n(M)\rOnto^{d_n} P_{n-1}(M).
$$
Let $[\omega_i]\in \pi_1(M\smallsetminus\{p_1,\ldots,p_{n-1}\})$ represented by a small circle around $p_i$. By Lemma~\ref{lemma3.10}, for each $1\leq i\leq n-1$, the subgroup
  $\Ker(d_i)\cap \Ker(d_n)$
is the normal closure  of $[\omega_i]$ in $\pi_1(M\smallsetminus\{p_1,\ldots,p_{n-1}\})$. Let $R_i=\Ker(d_i)\cap \Ker(d_n)$. Note that $\pi_1(M\smallsetminus\{p_1,\ldots,p_{n-1}\})$ is a free group with a basis containing the elements $[\omega_i]$ for $1\leq i\leq n-1$.
By Lemma~\ref{lemma4.2},
$$
\begin{array}{rcl}
\bigcap\limits_{i=1}^k R_i&=&[R_1,R_2,\ldots, R_k]_S\\
&=&\prod_{j=1}^k\left[\bigcap\limits_{i\in \{1,\dots, \hat j,\dots, k\}}R_i,R_j\right]\\
\end{array}
$$
because $R_i$ is the kernel of
$$
d_i|\colon \pi_1(M\smallsetminus \{p_1,\ldots,p_{n-1}\})\longrightarrow \pi_1(M\smallsetminus \{p_1,\ldots,p_{i-1},p_{i+1},\ldots,p_{n-1}\})
$$
for $1\leq i\leq n-1$, and
$$
\bigcap\limits_{i\in \{1,\dots, \hat j,\dots, k\}}R_i=[R_1,R_2,\ldots, R_{j-1},R_{j+1},\ldots,R_k]_S.
$$
It should be noticed also that for normal subgroups $H_1$, $H_2$
$H_3$ of a group $G$
$$
[H_1 , H_3] \, [ H_2, H_3] = [H_1 H_2, H_3] ,
$$
see, for example,  identity (2') of Proposition~1.1 in \cite{Serr}.
It follows that
\begin{equation}\label{equation5.1}
\begin{array}{rl}
&\bigcap_{i=1}^k(\Ker(d_i)\cap \Ker(d_n))\\
=&\prod_{j=1}^k\left[\bigcap\limits_{i\in\{1,\dots, \hat j,\dots, k\}}(\Ker (d_i)\cap \Ker(d_n)),\Ker(d_j)\cap \Ker(d_n)\right]\\
\leq& \prod_{j=1}^k\left[\bigcap\limits_{i\in \{1,\dots,\hat j,\dots, k\}}(\Ker (d_i)\cap \Ker(d_n)),\Ker(d_j)\right].\\
\end{array}
\end{equation}
On the other hand, since
$$
\left[\bigcap\limits_{i\in \{1,\dots,\hat j,\dots,k\}}(\Ker
(d_i)\cap \Ker(d_n)),\Ker(d_j)\right]\leq \bigcap_{i=1}^k
(\Ker(d_i)\cap \Ker(d_n)),
$$
for every $j=1,\dots, k$, we have
\begin{equation}\label{equation5.2}
\bigcap_{i=1}^k(\Ker(d_i)\cap
\Ker(d_n))=\prod_{j=1}^k\left[\bigcap\limits_{i\in\{1,\dots, \hat
j,\dots, k\}}(\Ker (d_i)\cap \Ker(d_n)),\Ker(d_j)\right].
\end{equation}
By induction, the subgroup $\bigcap\limits_{i\in\{1,\dots, \hat
j,\dots, k\}}(\Ker (d_i)\cap \Ker(d_n))$ has finitely many normal
generators for every $j=1,\dots, k$. From the short exact sequence
of groups
$$
\pi_1(M\smallsetminus\{p_1,\ldots,p_{k-1},p_{k+1},\ldots,p_{n-1}\})\rInto P_n(M)\rOnto^{d_k} P_{n-1}(M),
$$
the subgroup $\Ker(d_k)$ has finitely many generators. By Lemma~\ref{lemma5.1}, the commutator subgroup
$$\left[\bigcap\limits_{i\in\{1,\dots, \hat
j,\dots, k\}}(\Ker (d_i)\cap \Ker(d_n)),\Ker(d_j)\right]
$$
has finitely many normal generators for every $j=1,\dots, k$ and
hence the group $\bigcap_{i=1}^k(\Ker(d_i)\cap \Ker(d_n))$ has
finitely many normal generators. The induction is finished.
\end{proof}

\subsection{Proof of Theorem~\ref{theorem1.2}}
The proof is given by two different cases.

\noindent\textbf{Case 1.} $M$ is a connected compact manifold with nonempty boundary. By Lemma~\ref{lemma5.2},
$$
\Brun_n(M)=\bigcap_{i=1}^{n-1}\Ker(d_i)\cap \Ker(d_n)
$$
has finitely many normal generators in $P_n(M)$. Thus the factor groups $P_n(M)/\Brun_n(M)$ and $B_n(M)/\Brun_n(M)$ are finitely presented.

\noindent\textbf{Case 2.} $M$ is a compact closed manifold. Let $\tilde M=M\smallsetminus\{q_1\}$.
Using the exact sequence of the fibration of Theorem~\ref{theorem3.1} and induction on $n$ we conclude  that the inclusion $f\colon \tilde M\to M$ induces an epimorphism
$$
f^n_*\colon P_n(\tilde M)\to P_n(M).
$$
Since
$$
\Brun_n(\tilde M)=[\la\la A_{1,n}\ra\ra^{P_n(\tilde M)},\la\la
A_{2,n}\ra\ra^{P_n(\tilde M)},\ldots, \la\la
A_{n-1,n}\ra\ra^{P_n(\tilde M)}]_S,
$$
we have
$$
f^n_*(\Brun_n(\tilde M))=[\la\la A_{1,n}\ra\ra^{P_n(M)},\la\la
A_{2,n}\ra\ra^{P_n(M)},\ldots, \la\la A_{n-1,n}\ra\ra^{P_n(M)}]_S.
$$
From the fact that $\Brun_n(\tilde M)$ has finitely many normal generators in $P_n(\tilde M)$, the subgroup
$$[\la\la A_{1,n}\ra\ra^{P_n(M)},\la\la A_{2,n}\ra\ra^{P_n(M)},\ldots, \la\la
A_{n-1,n}\ra\ra^{P_n(M)}]_S
$$
has finitely many normal generators in $P_n(M)$.

If $M\not=S^2$ or $\RP^2$ with $n\geq 3$, then, by Theorem~\ref{theorem1.1} and Proposition~\ref{proposition3.6}, the subgroup
$$
\Brun_n(M)=[\la\la A_{1,n}\ra\ra^{P_n(M)},\la\la
A_{2,n}\ra\ra^{P_n(M)},\ldots, \la\la A_{n-1,n}\ra\ra^{P_n(M)}]_S
$$
has finitely many normal generators in $P_n(M)$. Thus $P_n(M)/\Brun_n(M)$ and $B_n(M)/\Brun_n(M)$ are finitely presented for $M\not=S^2$ or $\RP^2$ with $n\geq 3$.

If $M=S^2$, then $P_3(S^2)/\Brun_3(S^2)=\{1\}$ and $B_3(S^2)/\Brun_3(S^2)=\Z/2$. For $n=4$, $\Brun_4(S^2)$ has $5$ generators according to~\cite[Proposition 7.2.1]{BCWW}. Thus $P_4(S^2)/\Brun_4(S^2)$ and $B_4(S^2)/\Brun_4(S^2)$ are finitely presented. For $n\geq 5$, by Theorem~\ref{theorem1.1}, $\Brun_n(S^2)$ is a finite extension of the subgroup
$$
[\la\la A_{1,n}\ra\ra^{P_n(S^2)},\la\la
A_{2,n}\ra\ra^{P_n(S^2)},\ldots, \la\la
A_{n-1,n}\ra\ra^{P_n(S^2)}]_S
$$
because $\pi_{n-1}(S^2)$ is finite.  Thus $\Brun_n(S^2)$ has finitely many normal generators in $P_n(S^2)$ and so the assertion holds for the case  $M=S^2$.

If $M=\RP^2$, then $\Brun_3(\RP^2)$ has $9$ generators according to Proposition~\ref{brun3-rp2}. Thus $P_3(\RP^2)/\Brun_3(\RP^2)$ and $B_3(\RP^2)/\Brun_3(\RP^2)$ are finitely presented. For $n\geq 4$, by (3) of Theorem~\ref{theorem1.1} together with fact that $\pi_{n-1}(S^2)$ is finitely generated, the subgroup $\Brun_n(\RP^2)$ has finitely many normal generators and so the assertion holds for the case $M=\RP^2$.
\hfill $\Box$

\section{Some Remarks}\label{section6}
\subsection{An algorithm for determining a free basis for Brunnian Braids}
By Lemma~\ref{lemma3.10}, for having a free basis for $\Brun_{n+1}(M)$, it suffices to determine a free basis for
\begin{equation}\label{equation6.1}
\bigcap_{i=1}^{n}\Ker(d_i|\colon \pi_1(M\smallsetminus\{p_1,\ldots,p_n\})\to \pi_1(M\smallsetminus\{p_1,\stackrel{\wedge i}{\cdots\cdots},p_n\})).
\end{equation}
 Let $M$ be connected $2$-manifold with nonempty boundary and let $\omega_i$ be a small circle around $p_i$. Then
 $$
 \pi_1(M\smallsetminus\{p_1,\ldots,p_n\})=F(S_0\sqcup\{[\omega_1],\ldots,[\omega_n]\}),
 $$
where $\pi_1(M)=F(S_0)$. Let $S$ be a set and let $T$ be a subset of $S$. By a \textit{projection homomorphism}
$$
\pi\colon F(S)\to F(T)
$$
we mean here a group homomorphism such that
$$
\pi(x)=\left\{
\begin{array}{rcl}
x&\textrm{ if }& x\in T,\\
1&\textrm{ if }& x\not\in T.\\
\end{array}\right.
$$
Note that for each subset $T$ of $S$ there is a unique projection homomorphism $\pi\colon F(S)\to F(T)$. In our case, the homomorphisms $d_i|$ are projection homomorphisms in the following sense:

Let $S=S_0\sqcup\{[\omega_1],\ldots,[\omega_n]\}$ and let
$$T_i=S_0\sqcup\{[\omega_1],\ldots,[\omega_{i-1}], [\omega_{i+1}],\ldots, [\omega_n]\}$$
for $1\leq i\leq n$. Then
$$
d_i|\colon F(S)\to F(T_i)
$$
is the projection homomorphism for each $1\leq i\leq n$. The algorithm in ~\cite[Section 3]{Wu1} provides a recursive formula to determine a free basis for the intersection subgroup $\bigcap\limits_{i=1}^n\Ker(d_i|)$ given as follows:

For $x$ a reduced word in alphabet $S$  and $y$ a reduced word in
alphabet $T$, define $\mu(x,y)$ by induction on the word length of $y$:
\begin{enumerate}
\item[1)] $\mu(x,y)=x$ if $y$ is the empty word;
\item[2)] $\mu(x,y)=[\mu(x,y'),z^{\epsilon}]$ if $y=y'z^{\epsilon}$ with $z\in T$ and $\epsilon=\pm 1$.
\end{enumerate}

Let $V$ be a set of reduced words in alphabet $S$ and $W$ be
a set of reduced words in alphabet $T$, a sub-alphabet of $S$. Define a set of words in alphabet $S$:
$$
\mathcal{A}(V)_W=\{\mu(x,y) \ | \ x\in V\textrm{ and } y \in W\}.
$$
By~\cite[Proposition 3.3]{Wu1}, $\mathcal{A}(\{S\setminus T\})_{F(T)}$ is a free basis for the kernel of the projection homomorphism $\pi\colon F(S)\to F(T)$. Now for the subsets $T_1,\ldots,T_n$ of $S$, construct a subset
$$
\mathcal{A}(T_1,\ldots,T_k)
$$
of $F(S)$ by induction on $k$ for $1\leq k\leq n$:
\begin{enumerate}
\item $\mathcal{A}(T_1)=\mathcal{A}(\{S\setminus T_1\})_{F(T_1)}$.
\item Let
$$
T_2^{(2)}=\{w\in \mathcal{A}(T_1) \ | \ w=[[x,y_1^{\epsilon_1}],\ldots,y_t^{\epsilon_t}]\textrm{ with } x,y_j\in T_2\textrm{ for all } j\}
$$
and define
$$
\mathcal{A}(T_1,T_2)=\mathcal{A}(\mathcal{A}(T_1))_{F(T^{(2)}_2)}.
$$
\item Suppose that $\mathcal{A}(T_1,\ldots,T_{k-1})$ is defined such that all of the elements in $\mathcal{A}(T_1,\ldots,T_{k-1})$ are written as certain commutators in $F(S)$ in terms of elements in $S$. Let
$$
T^{(k)}_k=\{w\in \mathcal{A}(T_1,\ldots,T_{k-1}) \ | \ w=[x_1^{\epsilon_1},\ldots, x_l^{\epsilon_l}]\textrm{ with } x_j\in T_k \textrm{ for all } j\},
$$
where $[x^{\epsilon_1}_1,\ldots, x^{\epsilon_t}_t]$ are the elements in $\mathcal{A}(T_1,\ldots,T_{k-1})$ that are written as iterated commutators. Define
$$
\mathcal{A}(T_1,\ldots,T_k)=
\mathcal{A}(\mathcal{A}(T_1,\ldots,T_{k-1}))_{F(T^{(k)}_k)}.
$$
\end{enumerate}
By~\cite[Theorem 3.4]{Wu1}, $\mathcal{A}(T_1,\ldots,T_k)$ is a free basis for $\bigcap\limits_{i=1}^k\Ker(d_i|)$ for $1\leq k\leq n$. In particular,
$$
\mathcal{A}(T_1,\ldots,T_n)
$$
is a free basis for $\bigcap\limits_{i=1}^n \Ker(d_i|)$.

\subsection{Bi-$\Delta$-structure on braids, Brunnian braids and the inverse braid monoids}
Let $M$ be a connected manifold with $\partial M\not=\emptyset$. Let $a$ be a point in a collar of $\partial M$
$$
\partial M\times [0,1)\subseteq M.
$$
Then the map
$$
\begin{array}{c}
F(M,n)\simeq F(M\smallsetminus \partial M\times [0,1),n)\rTo F(M,n+1),\\
(z_1,\ldots,z_n)\mapsto (z_1,\ldots,z_{i-1},a,z_{i+1},\ldots,z_n)\\
\end{array}
$$
induces a group homomorphism
$$
d^i\colon B_n(M)=\pi_1(F(M,n)/\Sigma_n)\rTo B_{n+1}(M)=\pi_1(F(M,n+1)/\Sigma_{n+1})
$$
for $1\leq i\leq n$. Intuitively, $d^i$ is given by adding a trivial strand in position $i$. According to~\cite[Example 1.2.8]{Wu4}, the sequence of groups $\{B_{n+1}(M)\}_{n\geq 0}$ with faces relabeled as $\{d_0,d_1,\ldots\}$ and cofaces relabeled as $\{d^0,d^1,\ldots\}$ forms a bi-$\Delta$-set structure. Namely the following identities hold:
\begin{enumerate}
\item $d_jd_i=d_id_{j+1}$ for $j\geq i$;
\item $d^jd^i=d^{i+1}d^j$ for $j\leq i$;
\item $d_jd_i=\left\{
\begin{array}{lcl}
d^{i-1}d_j&\textrm{ if }& j<i,\\
\id&\textrm{ if }& j=i,\\
d^id_{j-1}&\textrm{ if }& j>i.\\
\end{array}
\right.$
\end{enumerate}
By restricting to pure braid groups, the sequence of pure braid groups $\{P_{n+1}(M)\}_{n\geq0}$ is a bi-$\Delta$-group. According to~\cite[Proposition 1.2.9]{Wu4}, we have the following decomposition.
\begin{prop}\label{proposition6.1}
Let $M$ be a connected $2$-manifold with nonempty boundary. Then $P_n(M)$ is the (iterated) semi-direct product of the subgroups
$$
d^{i_k}d^{i_{k-1}}\cdots d^{i_1}(\Brun_{n-k}(M)),
$$
$1\leq i_1<i_2<\cdots<i_k\leq n$, $0\leq k\leq n-1$, with lexicographical
order from the right.\hfill $\Box$
\end{prop}
The braids in the subgroup $d^{i_k}d^{i_{k-1}}\cdots d^{i_i}(\Brun_{n-k}(M))$ can be described as $(n-k)$-strand Brunnian braids with $k$ dots (or straight lines) in positions $i_1, \ldots, i_k$. This fits with terminology of inverse braid monoids~\cite{Vershinin}.

Let $M$ be any connected $2$-manifold. Define a set
$$
\mathfrak{H}^B_n(M)=\{\beta\in B_n(M)\ | \ d_1\beta=d_2\beta=\cdots=d_n\beta\}
$$
Namely $\mathfrak{H}^B_n(M)$ consists of $n$-strand pure braids such that it stays the same braid after removing any one of its strands.
A typical element in $\mathfrak{H}^B_n(M)$ is the half-twist braid
$$
\Delta_n=(\sigma_1\sigma_2\cdots\sigma_{n-1})(\sigma_1\sigma_2\cdots\sigma_{n-2})\cdots (\sigma_1\sigma_2)\sigma_1.
$$
\begin{prop}\label{proposition6.2}
Let $M$ be any connected $2$-manifold. Then the set
$\mathfrak{H}^B_n(M)$ is subgroup of $B_n(M)$. Moreover
$d_i(\mathfrak{H}^B_n(M))\subseteq \mathfrak{H}^B_{n-1}(M)$ and the function
$$
d_1=d_2=\cdots=d_n\colon \mathfrak{H}_n^B(M)\to \mathfrak{H}_{n-1}^B(M)
$$
is a group homomorphism.
\end{prop}
\begin{proof}
Let $\beta,\gamma\in \mathfrak{H}^B_n(M)$. Then
\begin{equation}\label{equation6.2}
d_i(\beta\gamma)=d_i(\beta) d_{i\cdot\beta}(\gamma)=d_1(\beta)d_1(\gamma)
=d_i(\beta)d_i(\gamma)
\end{equation}
for $1\leq i\leq n$. Thus $\beta\gamma\in \mathfrak{H}^B_n(M)$. From
$$
1=d_i(\beta^{-1}\beta)=d_i(\beta^{-1})d_{i\cdot\beta^{-1}}(\beta)=d_i(\beta^{-1})d_1(\beta),
$$
we have
$$
d_i(\beta^{-1})=(d_1(\beta))^{-1}
$$
for $1\leq i\leq n$. Thus $\mathfrak{H}^B_n(M)$ is a subgroup of $B_n(M)$.

Now let $\beta\in \mathfrak{H}^B_n(M)$. From the identity
$$
d_j(d_i\beta)=d_j(d_1\beta)=d_1(d_{j+1}\beta)=d_1(d_1\beta),
$$
we have $d_i\beta=d_1\beta\in \mathfrak{H}^B_{n-1}(M)$. From equation~(\ref{equation6.2}),
$$
d_1=d_i\colon \mathfrak{H}^B_n(M)\to \mathfrak{H}^B_{n-1}(M)
$$
is a group homomorphism and hence the result.
\end{proof}

\begin{prop}\label{proposition6.3}
Let $M$ be any connected $2$-manifold. Let $n\geq 2$. Then $\mathfrak{H}^B_n(M)\cap P_n(M)$ is a subgroup of $\mathfrak{H}^B_n(M)$ of index $2$.
\end{prop}
\begin{proof}
Consider the short exact sequence of groups
$$
P_n(M)\rInto B_n(M)\rOnto^{q_n}\Sigma_n.
$$
The face function $d_i\colon B_n(M)\to B_{n-1}(M)$ induces a unique face $$d_i^{\Sigma}\colon \Sigma_n\to \Sigma_{n-1}$$ such that
$$
d_i^{\Sigma}\circ q_n=q_{n-1}\circ d_i^B
$$
for each $1\leq i\leq n$. Let $\mathfrak{H}_n^{\Sigma}=q_n(\mathfrak{H}_n^B(M))$. Then there is a left exact sequence
$$
1\rTo \bigcap_{i=1}^n\Ker(d_i^{\Sigma})\rTo \mathfrak{H}^{\Sigma}_n\rTo^{d_1^{\Sigma}}\mathfrak{H}^{\Sigma}_{n-1}.
$$
By direct calculation, $\bigcap_{i=1}^n\Ker(d_i^{\Sigma})=\{1\}$ for $n\geq 3$. Thus
$$
d_1^{\Sigma}\circ \cdots \circ d_1^{\Sigma}\colon \mathfrak{H}^{\Sigma}_n\longrightarrow \mathfrak{H}^{\Sigma}_2=\Z/2
$$
is a monomorphism for $n\geq 3$. Since the half-twist braid $\Delta_n$ has nontrivial image in $\mathfrak{H}^{\Sigma}_n$, we have
$$
\mathfrak{H}^{\Sigma}_n=\Z/2
$$
for $n\geq 2$ and hence the result.
\end{proof}

Let $\mathfrak{H}_n(M)=\mathfrak{H}^B_n(M)\cap P_n(M)$. Then $d_1(\mathfrak{H}_n(M))\leq \mathfrak{H}_{n-1}(M)$. This gives a tower of groups
$$
\cdots \rTo^{d_1} \mathfrak{H}_n(M)\rTo^{d_1} \mathfrak{H}_{n-1}(M)\rTo^{d_1}\cdots.
$$
Let $\mathfrak{H}(M)=\lim\limits_n\mathfrak{H}_n(M)$ be the inverse limit of the tower of groups.
\begin{prop}\label{proposition6.4}
Let $M$ be any connected $2$-manifold such that $M\not=S^2$ or $\RP^2$. Then
$$
d_1\colon \mathfrak{H}_n(M)\to\mathfrak{H}_{n-1}(M)
$$
is an epimorphism for each $n\geq 2$.
\end{prop}
\begin{proof}
If $M$ is a connected $2$-manifold with nonempty boundary, the assertion follows from ~\cite[Proposition 1.2.1]{Wu4} by using the bi-$\Delta$-structure on $\{P_{n+1}(M)\}_{n\geq0}$.

Suppose that $M$ is a closed manifold with $M\not=S^2$ or $\RP^2$. Let $\tilde M=M\smallsetminus\{q_1\}$. We show by induction that
$$
\mathfrak{H}_k(\tilde M)\rTo \mathfrak{H}_k(M)
$$
is onto for each $k$. Clearly the statement holds for $k=1$. Assume that the statement holds for $k-1$ with $k\geq 2$. Consider the commutative diagram
\begin{diagram}
\Brun_{k}(\tilde M)&\rInto &\mathfrak{H}_k(\tilde M)&\rOnto^{d_1}&\mathfrak{H}_{k-1}(\tilde M)\\
\dTo&&\dTo&&\dOnto\\
\Brun_{k}(M)&\rInto &\mathfrak{H}_k(M)&\rTo^{d_1}&\mathfrak{H}_{k-1}(M),\\
\end{diagram}
where the right column is onto by induction. By Theorem~\ref{theorem1.1}, $$\Brun_k(\tilde M)\to \Brun_k(M)$$ is onto and so the middle column $\mathfrak{H}_k(\tilde M)\to \mathfrak{H}_k(M)$ is onto. The induction is finished and so
$$
\mathfrak{H}_n(\tilde M)\to \mathfrak{H}_{n}(M)
$$
is an epimorphism for each $n$. From the right square in the above diagram,
$$
d_1\colon \mathfrak{H}_n(M)\to \mathfrak{H}_{n-1}(M)
$$
is an epimorphism for each $n\geq 2$ and hence the result.
\end{proof}

\begin{cor}\label{corollary6.5}
Let $M$ be any connected $2$-manifold such that $M\not=S^2$ or $\RP^2$ and let $\alpha\in B_n(M)$. Then the equation
$$
d_1\beta=\cdots=d_{n+1}\beta=\alpha
$$
for $(n+1)$-strand braids $\beta$ has a solution if and only if $\alpha$ satisfies the condition that
$$
d_1\alpha=\ldots=d_n\alpha.
$$
\end{cor}
\begin{proof}
If there exists $\beta$ such that $d_1\beta=\cdots=d_{n+1}\beta=\alpha$. Then $\alpha\in \mathfrak{H}^B_n(M)$ by Proposition~\ref{proposition6.1} and so $d_1\alpha=\cdots=d_n\alpha$.

Conversely suppose that $d_1\alpha=\ldots=d_n\alpha$. Then $\alpha\in \mathfrak{H}^B_n(M)$. If $\alpha\in \mathfrak{H}_n(M)$, then the equation in the statement has a solution for $\alpha$ by Proposition~\ref{proposition6.4}. If $\alpha\not\in\mathfrak{H}_n(M)$, then $\Delta_n\alpha\in \mathfrak{H}_n(M)$. Thus there exists $\gamma$ such that $d_1\gamma=\cdots=d_{n+1}\gamma=\Delta_n\alpha$. Since $d_1\Delta_{n+1}=\cdots=d_{n+1}\Delta_{n+1}=\Delta_n$, we have
$$
d_1(\Delta_{n+1}^{-1}\gamma)=\cdots=d_{n+1}(\Delta_{n+1}^{-1}\gamma)=\alpha
$$
and hence the result.
\end{proof}

Let $M$ be a connected $2$-manifold with nonempty boundary. For $n\geq k$, the \textit{James-Hopf operation} is defined as a function $$H_{k,n}\colon \Brun_k(M)\rTo \mathfrak{H}_n(M)$$
by setting $H_{k,k}(\beta)=\beta$ with
\begin{equation}\label{equation6.4}
H_{k,n}(\beta)=\prod_{1\leq i_1<i_2<\cdots<i_{n-k}\leq n}d^{i_{n-k}}d^{i_{n-k-1}}\cdots d^{i_1}(\beta)
\end{equation}
with lexicographic order from right for $\beta\in \Brun_k(M)$. The function $H_{k,n}$ satisfies the property that
$$
d_iH_{k,n}(\beta) =H_{k,n-1}(\beta)
$$
for $1\leq i\leq n$, $k< n$ and $\beta\in \Brun_k(M)$. Thus the sequence of elements $\{H_{k,n}(\beta)\}$ lifts to $\mathfrak{H}(M)$ that defines a function
$$
H_{k,\infty}\colon \Brun_k(M)\rTo \mathfrak{H}(M).
$$
According to~\cite[Theorem 1.2.4]{Wu4}, we have the following proposition.

\begin{prop}\label{proposition6.5}
Let $M$ be a connected $2$-manifold with nonempty boundary and let $\alpha\in\mathfrak{H}_n(M)$ with $1\leq n\leq \infty$. Then there exists an unique element $\delta_k(\alpha)\in \Brun_k(M)$ for $1\leq k\leq n$ such that the equation
    $$
    \alpha=\prod_{k=1}^n H_{k,n}(\delta_k(\alpha))
    $$
holds.\hfill $\Box$
\end{prop}

\section{Acknowledgments}

This work was started during the visit of Jie Wu to the University Montpellier~2 in December 2003 and it was continued during the visit of
V.~V.~Vershinin to the National University of Singapore in July 2007.
The main results were obtained during the stay
of the four authors in Oberwolfach in the framework of the program RiP
during May 18th - June 7th, 2008. The authors would like to thank the Mathematical Institute of Oberwolfach for its hospitality.
The work was continued by the last two authors during the visit of the
 third author in the Institute for Mathematical Sciences of the National University of Singapore during the program on Algebraic Topology, Braids and Mapping Class Groups
December 4 - 19,  2008. He would like to express his gratitude to the
IMS.

The first three authors were partially supported by the Russian-Indian grant.

The last author is partially supported by the Academic Research Fund
of the National University of Singapore R-146-000-101-112

\end{document}